\documentclass[reqno,a4paper,11pt]{amsart}
\usepackage{graphicx}
\usepackage{geometry}
\usepackage{a4wide}

\usepackage{latexsym}
\usepackage{amsmath}
\usepackage{amsfonts}
\usepackage{amssymb}
\usepackage{amsthm}
\usepackage{mathrsfs}

\usepackage{xfrac}
\usepackage{tikz}
\usetikzlibrary{positioning,shapes,shadows,arrows}
\usepackage{mathtools}

\usepackage[pdftex]{hyperref}
\usepackage[utf8x]{inputenc}
\usepackage[english]{babel}

\newcommand{\jf}{\mathrm{J}^1}
\newcommand{\jfr}{\mathrm{J}_{fr}^1}
\newcommand{\triend}{\mbox{\hspace{0.2mm}}\hfill$\triangle$}
\newcommand{\HE}{\mathbb{E}\mathrm{xt}}

\newtheorem{theorem}{Theorem}[section]
\newtheorem{proposition}[theorem]{Proposition}

\newtheorem{corollary}[theorem]{Corollary}

\theoremstyle{remark}
\newtheorem{example}[theorem]{Example}

\theoremstyle{definition} 
\newtheorem{definition}[theorem]{Definition}

\theoremstyle{remark} 
\newtheorem{remark}[theorem]{Remark}

\title{Symplectic structures on moduli spaces of framed sheaves on surfaces}

\author{Francesco Sala}

\address{Department of Mathematics, School of Mathematical and Computer Sciences, Heriot-Watt University, Colin Maclaurin Building, Riccarton, Edinburgh EH14 4AS, United Kingdom}



\email{F.Sala@hw.ac.uk, salafra83@gmail.com}

\begin{document}

\begin{abstract}
We provide generalizations of the notions of Atiyah class and Kodaira-Spencer map to the case of framed sheaves. Moreover, we construct closed two-forms on the moduli spaces of framed sheaves on surfaces. As an application, we define a symplectic structure on the moduli spaces of framed sheaves on some birationally ruled surfaces.
\end{abstract}

\maketitle

\setlength{\parskip}{0.2ex} 

\tableofcontents

\setlength{\parskip}{0.7ex plus 0.4ex minus 0.2ex} 

\section{Introduction}

Let $\mathcal{M}(r,n)$ be the moduli space of \textit{framed sheaves} on $\mathbb{CP}^2$, that is, the moduli space of pairs $(E,\alpha)$ modulo isomorphism, where $E$ is a torsion-free sheaf on $\mathbb{CP}^2$ of rank $r$ with $c_2(E)=n$, locally trivial in a neighborhood of a fixed line $l_{\infty}$, and $\alpha\colon E\vert_{l_{\infty}}\stackrel{\sim}{\rightarrow} \mathcal{O}_{l_{\infty}}^{\oplus r}$ is the \textit{framing at infinity}. $\mathcal{M}(r,n)$ is a nonsingular quasi-projective variety of dimension $2rn.$ This moduli space also admits a description in terms of linear data, the so-called \textit{ADHM} data (see, for example, Chapter 2 of Nakajima's book \cite{book:nakajima1999}). 

As described in \cite[Chapter 3]{book:nakajima1999}, by using the ADHM data description, the moduli space $\mathcal{M}(r,n)$ can be realized as a hyper-K\"ahler quotient. By fixing a complex structure within the hyper-K\"ahler family of complex structures on $\mathcal{M}(r,n)$, one can define a holomorphic symplectic form on $\mathcal{M}(r,n).$ 

In addition to the case of $\mathbb{CP}^2$, the only other relevant result in the literature about symplectic structures on moduli spaces of framed sheaves is due to Bottacin \cite{art:bottacin2000}. 

Let $X$ be a complex nonsingular projective surface, $D$ an effective divisor such that $D=\sum_{i=1}^n C_i$, where $C_i$ is an integral curve for $i=1,\ldots, n$, and $F_D$ a locally free $\mathcal{O}_D$-module. We call $(D,F_D)$-\emph{framed vector bundle} on $X$ a pair $(E,\alpha)$, where $E$ is a locally free sheaf on $X$ and $\alpha\colon E\vert_D\stackrel{\sim}{\rightarrow} F_D$ is an isomorphism. Let us fix a Hilbert polynomial $P.$ Bottacin constructs Poisson structures on the moduli space $\mathcal{M}_{lf}^*(X;F_D,P)$ of framed vector bundles on $X$ with Hilbert polynomial $P$, that are induced by global sections of the line bundle $\omega_X^{-1}(-2\,D).$ In particular, when $X$ is the complex projective plane, $D=l_\infty$ is a line and $F_D$ the trivial vector bundle of rank $r$ on $l_\infty$, this yields a symplectic structure on the moduli space $\mathcal{M}_{lf}(r,n)$ of framed vector bundles on $\mathbb{CP}^2$, induced by the standard holomorphic symplectic structure of $\mathbb{C}^2=\mathbb{CP}^2\setminus l_\infty.$ It is not known if this symplectic structure is equivalent to the one given by the ADHM construction.

Bottacin's result can be seen as a generalization to the framed case of the construction of Poisson brackets and holomorphic symplectic two-forms on the moduli spaces of Gieseker-stable torsion-free sheaves on $X.$ We recall briefly the main results for torsion-free sheaves. In \cite{art:mukai1984}, Mukai proved that any moduli space of simple sheaves on a $K3$ or abelian surface has a non-degenerate holomorphic two-form. Its closedness was later proved by Mukai in \cite{art:mukai1987}. Mukai's result was generalizated to moduli spaces of simple vector bundles on symplectic K\"ahler manifolds by Ran \cite{art:ran1996} and to moduli spaces of Gieseker-stable vector bundles over surfaces of general type and over Poisson surfaces by Tyurin \cite{art:tyurin1988}; a more thorough study of the Poisson case was made by Bottacin in \cite{art:bottacin1995}. In \cite{art:ogrady1992}, by using these results O' Grady defined closed two-forms on algebraic varieties parametrizing flat families of coherent sheaves. In all these cases, the symplectic two-form is defined in terms of the Atiyah class.

In the present paper we define a modified Atiyah class of a family of framed sheaves, which allows us to describe a framed version of the Kodaira-Spencer map and to construct closed two-forms on the moduli spaces of framed sheaves, that under some conditions are symplectic.

More precisely, let $X$ be a nonsingular projective surface over an algebraically closed field $k$ of characteristic zero, $D\subset X$ a divisor, $F_D$ a locally free $\mathcal{O}_D$-module and $S$ a Noetherian $k$-scheme of finite type. A flat family of $(D,F_D)$-framed sheaves parametrized by $S$ is a pair $\mathcal{E}:=(E,\alpha)$ on $S\times X$, such that $E$ is flat over $S$ and all the restrictions to the fibres $\{s\}\times X$ are $(D,F_D)$-framed sheaves on $X.$

Let $\mathcal{E}=(E,\alpha)$ be a $S$-flat family of $(D,F_D)$-framed sheaves. We introduce the \emph{framed sheaf of first jets} $\jfr(\mathcal{E})$ as the subsheaf of the sheaf of first jets $\jf(E)$ (introduced by Atiyah in \cite{art:atiyah1957}) consisting of those sections whose $p_S^*(\Omega_S^1)$-part vanishes along $S\times D.$ We define the \emph{framed Atiyah class} $at(\mathcal{E})$ of $\mathcal{E}$ as an extension class of $\jfr(\mathcal{E})$ in 
\begin{equation*}
\mathrm{Ext}^1(E,\left(p_S^*(\Omega^1_S)(-S\times D)\oplus p_X^*(\Omega_X^1)\right)\otimes E),
\end{equation*}
where $p_S^*(\Omega^1_S)(-S\times D)=p_S^*(\Omega^1_S)\otimes\mathcal{O}_{S\times X}(-S\times D)=p_S^*(\Omega^1_S)\otimes p_X^*(\mathcal{O}_X(-D)).$

Starting from the framed Atiyah class $at(\mathcal{E})$, one can define a section $\mathcal{A}t_S(\mathcal{E})$ in 
\begin{equation*}
\mathrm{H}^0(S,\mathcal{E}xt_{p_S}^1(E, p_S^*(\Omega^1_S)\otimes p_X^*(\mathcal{O}_X(-D))\otimes E)),
\end{equation*}
where $p_S\colon S\times X \rightarrow S$ is the projection.

In the same way as in the nonframed case (cf. \cite[Section 10.1.8]{book:huybrechtslehn2010}), by using  $\mathcal{A}t_S(\mathcal{E})$ one can define the \emph{framed} Kodaira-Spencer map associated to $\mathcal{E}$:
\begin{equation*}
KS_{fr}\colon (\Omega_S^1)^\vee\longrightarrow \mathcal{E}xt_{p_S}^1(E,p_X^*(\mathcal{O}_X(-D))\otimes E).
\end{equation*}
$(D,F_D)$-framed sheaves are a particular case of \emph{framed modules}, whose theory was developed by Huybrechts and Lehn \cite{art:huybrechtslehn1995-I, art:huybrechtslehn1995-II}. In order to construct moduli spaces parametrizing these objects, they defined a notion of semistability depending on a polarization and a rational polynomial $\delta$ of degree less than the dimension of the ambient variety and positive leading coefficient. They proved that there exists a coarse moduli space parametrizing semistable framed modules and a fine moduli space parametrizing stables ones. Moduli spaces of stable $(D,F_D)$-framed sheaves turn out to be open subschemes of the fine moduli spaces of stable framed modules (cf. \cite{art:bruzzomarkushevich2011}).

Let $X$ be a nonsingular projective surface over $k$ equipped with an ample line bundle $\mathcal{O}_X(1).$ We consider as above pairs $(D,F_D).$ Let $\delta\in \mathbb{Q}[n]$ be a stability polynomial and $P$ a numerical polynomial of degree two. 

Let us denote by $\mathcal{M}^*_\delta(X;F_D,P)$ the moduli space of $(D,F_D)$-framed sheaves on $X$ with Hilbert polynomial $P$ that are stable with respect to $\mathcal{O}_X(1)$ and $\delta.$ Let $\mathcal{M}^*_\delta(X;F_D,P)^{sm}$ be the smooth locus of $\mathcal{M}^*_\delta(X;F_D,P)$ and  $\tilde{\mathcal{E}}=(\tilde{E},\tilde{\alpha})$ a \emph{universal object} over $\mathcal{M}^*_\delta(X;F_D,P)^{sm}.$ 

The first result we obtain by using the framed Atiyah class is the following:
\begin{theorem}
The framed Kodaira-Spencer map defined by $\tilde{\mathcal{E}}$ induces a canonical isomorphism
\begin{equation*}
KS_{fr}\colon T\mathcal{M}^*_\delta(X;F_D,P)^{sm}\stackrel{\sim}{\longrightarrow} \mathcal{E}xt_p^1(\tilde{E}, \tilde{E}\otimes p_X^*(\mathcal{O}_X(-D))),
\end{equation*}
where $p$ is the projection from $\mathcal{M}^*_\delta(X;F_D,P)^{sm}\times X$ to $\mathcal{M}^*_\delta(X;F_D,P)^{sm}.$
\end{theorem}
Thus we get a generalization to the framed case of the corresponding statement for the moduli space of Gieseker-stable torsion-free sheaves on $X$ (cf. \cite[Theorem 10.2.1]{book:huybrechtslehn2010}).

From this theorem it follows that for any point $[(E,\alpha)]$ of $\mathcal{M}^*_\delta(X;F_D,P)^{sm}$, the vector space $\mathrm{Ext^1}(E,E(-D))$ is naturally identified with the tangent space $T_{[(E,\alpha)]}\mathcal{M}^*_\delta(X;F_D,P).$ For any $\omega\in \mathrm{H}^0(X,\omega_X(2\,D))$, we can define a skew-symmetric bilinear form
\begin{eqnarray*}
&&\mathrm{Ext}^1(E,E(-D))\times \mathrm{Ext}^1(E,E(-D))\stackrel{\circ}{\longrightarrow} \mathrm{Ext}^2(E,E(-2\,D))\\
&&\stackrel{tr}{\longrightarrow} \mathrm{H}^2(X, \mathcal{O}_X(-2\,D))\stackrel{\cdot\, \omega}{\longrightarrow} \mathrm{H}^2(X,\omega_X)\cong k.
\end{eqnarray*}
By varying $[(E,\alpha)]$, these forms define an exterior two-form $\tau(\omega)$ on $\mathcal{M}^*_\delta(X;F_D,P)^{sm}.$ 

We prove that $\tau(\omega)$ is a closed form (cf. Theorem \ref{thm:closed}) and provide a criterion of its non-degeneracy (cf. Proposition \ref{prop:non-dege}). In particular, if the line bundle $\omega_X(2\,D)$ is trivial, the two-form $\tau(1)$ induced by $1\in \mathrm{H}^0(X,\omega_X(2\,D))\cong k$ defines a symplectic structure on $\mathcal{M}^*_\delta(X;F_D,P)^{sm}.$ 

As an application, we construct holomorphic symplectic structures on moduli spaces of $(D,F_D)$-framed sheaves on some birationally ruled surfaces. In particular, for $\mathbb{CP}^1\times \mathbb{CP}^1$, $\mathbb{C}\times \mathbb{CP}^1$, with $C$ elliptic curve, and the second Hirzebruch surface $\mathbb{F}_2$, under a suitable choice of $D$ and for $F_D=\mathcal{O}_D^{\oplus r}$, we get a generalization to the non-locally free case of Bottacin's construction of symplectic structures induced by non-degenerate Poisson structures (cf. \cite{art:bottacin2000}).

For $X=\mathbb{CP}^2$ or the blowup of $\mathbb{CP}^2$ at a point, under a suitable choice of $D$ and $F_D$, we obtain new examples of non-compact holomorphic symplectic varieties, not covered by Bottacin's construction.

This paper is structured as follows. In Section \ref{sec:preliminaries}, we briefly introduce the theory of framed modules and framed sheaves and state the main theorems about the existence of moduli spaces parametrizing these objects. In Section \ref{sec:atiyah} we recall the definition of the Atiyah class for a flat family of coherent sheaves. In Section \ref{sec:framedatiyah} we give the definitions of the framed version of the Atiyah class and of the Kodaira-Spencer map. In Section \ref{sec:framedkodaira} we prove that the framed Kodaira-Spencer map is an isomorphism for the moduli space of stable $(D,F_D)$-framed sheaves and, in Section \ref{sec:framedtwoforms}, we construct closed two-forms on it. Finally, in Section \ref{sec:example} we apply our result to some particular birationally ruled surfaces, and define a symplectic structure on the moduli spaces of (stable) $(D,F_D)$-framed sheaves on them.

\subsection*{Conventions}

All schemes are Noetherian of finite type over an algebraically closed field $k$ of characteristic zero. A \emph{variety} is a reduced separated scheme. A \emph{polarized variety of dimension} $d$ is a pair $(X,\mathcal{O}_X(1))$, where $X$ is a smooth connected projective variety of dimension $d$, defined over $k$, and $\mathcal{O}_X(1)$ a very ample line bundle. The \emph{canonical line bundle} of $X$ is denoted by $\omega_X$ and its associated divisor class by $K_X.$

Let $(X,\mathcal{O}_X(1))$ be a polarized variety and $S$ a scheme. We denote by $\mathcal{X}$ the cartesian product $S\times X$, and by $p_S$, $p_X$ the projections from $\mathcal{X}$ to $S$ and $X$ respectively. Let $D$ be an effective divisor on $X$: we denote by $\mathcal{D}$ the cartesian product $S\times D$ and for any coherent sheaf $F$ on $\mathcal{X}$, $F(-\mathcal{D})$ is the tensor product $F\otimes p_X^*(\mathcal{O}_X(-D)).$

We denote by $NS(X)$ the Néron-Severi group of a smooth connected projective variety $X$, that is, the image of the map $c_1\colon \mathrm{Pic}(X)\rightarrow \mathrm{H}^2(X,\mathbb{Z}).$

As usual, a polarized variety of dimension two is called a \emph{polarized surface}.

Finally, we denote by $\HE^i(E^\bullet, G^\bullet)$ the $i$-th hyper-Ext group of two finite complexes of locally free sheaves $E^\bullet$ and $G^\bullet$ on a scheme $Y$, that is, the $i$-th hyper-cohomology group of the total complex associated to the double complex $C^\bullet(\mathcal{H}om^\bullet(E^\bullet, G^\bullet)).$ As it is explained in \cite[Section 10.1]{book:huybrechtslehn2010}, for any coherent sheaf $N$ on $Y$, one can define the trace map $\mathrm{tr}\colon \HE^i(E^\bullet,E^\bullet\otimes N)\rightarrow \mathrm{H}^i(Y,N).$

Moreover, if $E$ is a coherent sheaf on $Y$ that admits a finite locally free resolution $E^\bullet\rightarrow E$, we have $\HE^i(E^\bullet, E^\bullet\otimes N)\cong \mathrm{Ext}^i(E,E\otimes N)$ for any locally free sheaf $N.$ In this case, we denote by $\mathrm{Ext}^i(E,E\otimes N)_0$ the kernel
\begin{equation*}
\ker[\mathrm{tr}\colon \mathrm{Ext}^i(E,E\otimes N)\rightarrow \mathrm{H}^i(Y,N)].
\end{equation*}

\subsubsection*{Acknowledgements}

This paper is largely based upon the author’s PhD thesis \cite{phd:sala2011}. The author thanks his supervisors Ugo Bruzzo and Dimitri Markushevich for suggesting this problem and for their constant support. He thanks the referees for useful remarks and suggestions which helped to improve the paper. He also thanks Francesco Bottacin, Laurent Manivel and Christoph Sorger for interesting discussions. This paper was mostly written while the author was staying at SISSA, Université Lille 1 and IHÉS. He thanks those institutions for hospitality and support. He was partially supported by PRIN ``Geometria delle varietà algebriche e dei loro spazi di moduli'', co-funded by MIUR (cofin 2008), by the grant of the French Agence Nationale de Recherche VHSMOD-2009 Nr. ANR-09-BLAN-0104, by the European Research Network ``GREFI-GRIFGA'' and by the ERASMUS ``Student Mobility for Placements'' Programme.

\section{Preliminaries on framed sheaves}\label{sec:preliminaries}

In this section we introduce the notions of $(D,F_D)$-framed sheaves and framed modules. Moreover, we recall the main results about the construction of moduli spaces of these objects.

Let $(X,\mathcal{O}_X(1))$ be a polarized surface.
\begin{definition}\label{def:locfree}
Let $D$ be an effective divisor of $X$ and $F_D$ a coherent sheaf on $X$, which is a locally free $\mathcal{O}_{D}$-module. We say that a coherent sheaf $E$ on $X$ is $(D,F_D)$\textit{-framable} if $E$ is torsion-free, locally free in a neighborhood of $D$, and there is an isomorphism $E\vert_D\stackrel{\sim}{\rightarrow} F_D.$ An isomorphism $\alpha\colon E\vert_D\stackrel{\sim}{\rightarrow} F_D$ will be called a $(D,F_D)$-\textit{framing} of $E.$ A $(D,F_D)$-\textit{framed sheaf} is a pair $\mathcal{E}:=(E,\alpha)$ consisting of a $(D,F_D)$-framable sheaf $E$ and a $(D,F_D)$-framing $\alpha.$ Two $(D,F_D)$-framed sheaves $(E,\alpha)$ and $(E',\alpha')$ are isomorphic if there is an isomorphism $f\colon E\rightarrow E'$ such that $\alpha'\circ f\vert_D=\alpha.$
\end{definition}
\begin{remark}
Our notion of framing is the same as the one used in M. Lehn's Ph.D. thesis \cite{phd:lehn1993} and in T. A. Nevins' papers \cite{art:nevins2002-I, art:nevins2002-II}. \triend 
\end{remark}
We shall call $(D,F_D)$\emph{-framed vector bundles} on $X$ those $(D,F_D)$-framed sheaves whose underlying coherent sheaf is locally free.

Note that $(D,F_D)$-framed sheaves are a particular type of \emph{framed modules}.
\begin{definition}[\cite{art:huybrechtslehn1995-I, art:huybrechtslehn1995-II}]
Let $F$ be a coherent sheaf on $X.$ A $F$-\textit{framed module} on $X$ is a pair $(E,\alpha)$, where $E$ is a coherent sheaf on $X$ and $\alpha\colon E\rightarrow F$ a morphism of coherent sheaves.
\end{definition}
In \cite{art:huybrechtslehn1995-I, art:huybrechtslehn1995-II}, Huybrechts and Lehn developed the theory of framed modules. In the case of polarized surfaces, they introduced a notion of (semi)stability depending on the polarization $\mathcal{O}_X(1)$ and a rational polynomial $\delta(n)=\delta_1 n +\delta_0$, with $\delta_1>0.$ We shall call $\delta$ a \emph{stability polynomial}. When the framing $\alpha$ is zero, this reduces to Gieseker's (semi)stability condition for torsion-free sheaves.

Fix a coherent sheaf $F$ and a numerical polynomial $P$ of degree two. Let us denote by $\underline{\mathcal{M}}^{ss}_\delta(X;F,P)$ (resp. $\underline{\mathcal{M}}^{s}_\delta(X;F,P)$) the contravariant functor from the category of schemes to the category of sets, that associates with every scheme $S$ the set of isomorphism classes of flat families $(G,\beta\colon G\rightarrow p_X^*(F))$ of semistable (resp. stable) $F$-framed modules with Hilbert polynomial $P$ parametrized by $S.$ The main result in their papers is the following:
\begin{theorem}[\cite{art:huybrechtslehn1995-I, art:huybrechtslehn1995-II}]
There exists a projective scheme $\mathcal{M}^{ss}_\delta(X;F,P)$ which corepresents the moduli functor $\underline{\mathcal{M}}^{ss}_\delta(X;F,P)$ of semistable $F$-framed modules. Moreover, there is an open subscheme $\mathcal{M}^{s}_\delta(X;F,P)$ of $\mathcal{M}^{ss}_\delta(X;F,P)$ that represents the moduli functor $\underline{\mathcal{M}}^{s}_\delta(X;F,P)$ of stable $F$-framed modules, i.e., $\mathcal{M}^{s}_\delta(X;F,P)$ is a fine moduli space for stable $F$-framed modules.
\end{theorem}
``Fine'' means the existence of a \emph{universal} $F$-framed module over $\mathcal{M}^{s}_\delta(X;F,P)$, that is, a pair $(\tilde{E},\tilde{\alpha}\colon \tilde{E}\rightarrow p_X^{*}(F))$, where $\tilde{E}$ is a coherent sheaf on $\mathcal{M}^{s}_\delta(X;F,P)\times X$, flat over  $\mathcal{M}^{s}_\delta(X;F,P)$, such that for any scheme $S$ and any family of stable $F$-framed modules $(G,\beta)\in \underline{\mathcal{M}}^{s}_\delta(X;F,P)(S)$ parametrized by $S$, there exists a unique morphism $g\colon S\rightarrow \mathcal{M}^{s}_\delta(X;F,P)$ such that the pull back of $(\tilde{E},\tilde{\alpha})$ is isomorphic to $(G,\beta).$

By using Huybrechts and Lehn's result, Bruzzo and Markushevich constructed a moduli space parametrizing \emph{all} isomorphism classes of $(D,F_D)$-framed sheaves on $X$, under some mild assumptions on the divisor $D$ and on $F_D.$ More precisely, they proved the following:
\begin{theorem}[\cite{art:bruzzomarkushevich2011}]\label{thm:bruzzomarkushevich}
Let $X$ be a surface, $D\subset X$ a big and nef curve and $F_D$ a Gieseker-semistable locally free $\mathcal{O}_D$-module. Then for any numerical polynomial $P$ of degree two, there exists a polarization $\mathcal{O}_X(1)$ and a stability polynomial $\delta(n)=\delta_1 n+\delta_0$ for which there is an open subset $\mathcal{M}^*(X;F_D,P)\subset \mathcal{M}^{s}_\delta(X;F_D,P)$ of the moduli space of stable $F_D$-framed modules, which is a fine moduli space for $(D,F_D)$-framed sheaves on $X$ with Hilbert polynomial $P.$ In particular, $\mathcal{M}^*(X;F_D,P)$ is a quasi-projective separated scheme.
\end{theorem}

Now we give a characterization of the smoothness of the moduli space $\mathcal{M}^*(X;F_D,P).$
\begin{theorem}[\cite{art:bruzzomarkushevich2011}]
Let $[(E,\alpha)]$ be a point in $\mathcal{M}^*(X;F_D,P).$ The following statements hold:
\begin{itemize}
 \item[(i)] The Zariski tangent space of $\mathcal{M}^*(X;F_D,P)$ at $[(E,\alpha)]$ is naturally isomorphic to the Ext-group $\mathrm{Ext}^1(E,E(-D)).$
\item[(ii)] If $\mathrm{Ext}^2(E,E(-D))_0$ vanishes, then $\mathcal{M}^*(X;F_D,P)$ is smooth at $[(E,\alpha)].$
\end{itemize}
\end{theorem}
\begin{corollary}\label{cor:smoothness}
Let $X$ be a rational surface over $\mathbb{C}$, $D\subset X$ a smooth connected big and nef curve of genus zero and $F_D\cong \mathcal{O}_D^{\oplus r}.$ Let $n\in \mathbb{Z}$ and $c\in NS(X)$ such that $\int_{D} c=0.$ Then the moduli space $\mathcal{M}^*(X;F_D,r,c,n)$ parametrizing isomorphism classes of $(D,F_D)$-framed sheaves on $X$ of rank $r$, first Chern class $c$ and second Chern class $n$ is a smooth quasi-projective variety of dimension
\begin{equation*}
\dim_{\mathbb{C}} \mathcal{M}^*(X;F_D,r,c,n)=2rn-(r-1)\int_X c^2.
\end{equation*}
\end{corollary}
\proof
By \cite[Proposition 2.1]{art:gasparimliu2010}, for any point $[(E,\alpha)]\in \mathcal{M}^*(X;F_D,P)$ we have
\begin{equation*}
\mathrm{Ext}^0(E,E(-D))=\mathrm{Ext}^2(E,E(-D))=0.
\end{equation*}
Hence the moduli space $\mathcal{M}^*(X;F_D,P)$ is smooth. The formula for the dimension follows from the Hirzebruch-Riemann-Roch theorem.
\endproof

\section{The Atiyah class}\label{sec:atiyah}

In this section we recall the notion of the \emph{Atiyah class} for flat families of coherent sheaves. The Atiyah class was introduced by Atiyah in \cite{art:atiyah1957} for the case of coherent sheaves and by Illusie in \cite{book:illusie1971, book:illusie1972} for any complex of coherent sheaves (see \cite[Section 10.1.5]{book:huybrechtslehn2010} for a description of the Atiyah class in terms of \v{C}ech cocycles). Atiyah's approach involves the notion of the \emph{sheaf of first jets} of a fixed coherent sheaf associated to the sheaf of one-forms (for a generalization of the sheaf of first jets to quotients of the sheaf of one-forms, see Maakestad's paper \cite{art:maakestad2010-I}).

Let $Y$ be a scheme.
\begin{definition}
Let $E$ be a coherent sheaf on $Y.$ We call \textit{sheaf of the first jets} $\jf(E)$ of $E$ the coherent sheaf of $\mathcal{O}_Y$-modules defined as follows:
\begin{itemize}
 \item as a sheaf of $k$-modules, we set $\jf(E):=(\Omega_Y^1\otimes E)\oplus E$,
\item for any $y\in Y$, $a\in \mathcal{O}_{Y,y}$ and $(z\otimes e, f)\in \jf(E)_y$, we define 
\begin{equation*}
a(z\otimes e, f):=(az\otimes e +d(a)\otimes f, af),
\end{equation*}
where $d$ is the exterior differential of $Y.$
\end{itemize}
\end{definition}
The sheaf $\jf(E)$ fits into an exact sequence of coherent sheaves 
\begin{equation}\label{eq:atiyahclass}
0\longrightarrow \Omega^1_Y\otimes E \longrightarrow \jf(E)\longrightarrow E\longrightarrow 0.
\end{equation}
\begin{definition}
Let $E$ be a coherent sheaf on $Y.$ We call \textit{Atiyah class} of $E$ the class $at(E)$ in $\mathrm{Ext}^1(E,\Omega_Y^1\otimes E)$ associated to the extension \eqref{eq:atiyahclass}.
\end{definition}
It is a well known fact that the Atiyah class $at(E)$ is the obstruction for the existence of an algebraic connection on $E.$ This means the following:
\begin{proposition}{\normalfont (\cite[Proposition 3.4]{art:maakestad2010-I}).}
Let $E$ be a coherent sheaf on $Y.$ The Atiyah class $at(E)$ is zero if and only if there exists a connection on $E.$ 
\end{proposition}

\section{The Atiyah class for framed sheaves}\label{sec:framedatiyah}

In this section we turn to the framed case. First, we introduce the \emph{framed sheaf of first jets} $\jfr(\mathcal{E})$ of a $S$-flat family $\mathcal{E}=(E,\alpha)$ of $(D,F_D)$-framed sheaves as the subsheaf of the sheaf of first jets $\jf(E)$ consisting of those sections whose $p_S^*(\Omega_S^1)$-part vanishes along $S\times D.$ Then we define a \emph{framed} version of the Atiyah class and of the Newton polynomials\footnote{In the nonframed case, the Newton polynomials are introduced, for example, in \cite[Section 10.1.6]{book:huybrechtslehn2010}.} for $S$-flat families. We use this \emph{relative} approach, because later on we shall want to consider the framed Atiyah class of a \emph{universal framed sheaf}.

From now on we fix the pair $(D,F_D)$ and we just say \emph{framed sheaf} for a $(D,F_D)$-framed sheaf.
\begin{definition}
Let $S$ be a scheme. A \textit{flat family of framed sheaves parametrized by} $S$ is a pair $\mathcal{E}=(E,\alpha)$ where $E$ is a coherent sheaf on $\mathcal{X}$, flat over $S$, and $\alpha\colon E\rightarrow p_X^*(F_D)$ is a morphism such that for any $s\in S$ the pair $(E\vert_{\{s\}\times X},\alpha\vert_{\{s\}\times X})$ is a $\left(\{s\}\times D,p_X^*(F_D)\vert_{\{s\}\times D}\right)$-framed sheaf on $\{s\}\times X.$ 
\end{definition}
Let $\mathcal{E}=(E,\alpha)$ be a flat family of framed sheaves parametrized by a scheme $S.$ We define a subsheaf $\jfr(\mathcal{E})$ of $\jf(E)$, that we shall call \emph{framed sheaf of first jets} of $\mathcal{E}.$

For a point $x\in \mathcal{X}$ such that $x\notin \mathcal{D}$, we set $\jfr(\mathcal{E})_x:=\jf(E)_x.$ 

Fix $x\in \mathcal{D}.$ By definition of a flat family of framed sheaves, $\left(E\vert_{\{s\}\times X}\right)_x$ is a free module for $s\in S$ such that $p_S(x)=s$, hence $E_x$ is a free module (cf. \cite[Lemma 2.1.7]{book:huybrechtslehn2010}). Therefore there exists an open neighborhood $V\subset \mathcal{X}$ of $x$ such that $E\vert_V$ is a locally free $\mathcal{O}_V$-module. 

We denote by $E_V$ the restriction of $E$ to $V.$ Let $\mathcal{D}':=V\cap \mathcal{D}$ and $\mathcal{U}=\{U_i\}_{i\in I}$ a cover of $\mathcal{D}'$ over which $p_X^{*}(F_D)\vert_{\mathcal{D}'}$ trivializes, and choose on any $U_i$ a set $\{e_i^0\}$ of basis sections of $\Gamma(p_X^{*}(F_D)\vert_{\mathcal{D}'},U_i).$ Let $g^0_{ij}$ be transition functions of $p_X^{*}(F_D)\vert_{\mathcal{D}'}$ with respect to chosen local basis sections (i.e., $e^0_i=g^0_{ij}e^0_j$), constant along $S.$ Let us fix a cover $\mathcal{W}=\{W_i\}_{i\in I}$ of $V$ over which $E_V$ trivializes with sets $\{e_i\}$ of basis sections such that $W_i\cap \mathcal{D}'=U_i$ for any $i\in I$ and
\begin{eqnarray*}
e_i\vert_{\mathcal{D}'}&=& e_i^0, \\
g_{ij}\vert_{\mathcal{D}'}&=& g_{ij}^0.
\end{eqnarray*}
Let $z^1_i, \ldots, z^s_i$ and $z^{s+1}_i, \ldots, z^t_i$ be the local coordinates of $S$ and $X$ on $W_i$, respectively, and $f_i=0$ the local equation of $\mathcal{D}$ on $W_i.$ Define $\jfr(\mathcal{E})_x\subset \jf(E)_x=(\Omega_{\mathcal{X},x}^1\otimes E_x)\oplus E_x$ as the $\mathcal{O}_{\mathcal{X},x}$-module spanned by the basis obtained by tensoring all the elements of the set $\{f_i d z^1_i, \ldots, f_i d z^s_i, d z^{s+1}_i,\ldots, d z^t_i\}$, where $\{dz^1_i, \ldots, dz^t_i\}$ is a basis of $\Omega_{\mathcal{X},x}^1$, by the elements of the basis $\{e_i\}:=\{e^1_i, \ldots, e^r_i\}$ of $E_x$ and then adding the elements of $\{e_i\}.$ Thus an arbitrary element of $\jfr(\mathcal{E})_x$ is of the form
\begin{equation*}
h_i + f_i \sum_{n=1}^r \sum_{k=1}^s \chi_{n,k}\, e_i^n\otimes d z_i^k+\sum_{m=1}^r \sum_{l=s+1}^t \psi_{m,l}\, e_i^m\otimes d z_i^l,
\end{equation*}
where $h_i\in E_x$ and $\chi_{n,k}, \psi_{m,l}\in \mathcal{O}_{X,x}$, for $m,n=1, \ldots r$, $k=1,\ldots, s$ and $l=s+1,\ldots, t.$

If $x$ is also a point in the open subset $W_j$ of $\mathcal{W}$, let us denote by $l_{ij}\in \mathcal{O}_V^*(W_i\cap W_j)$ the transition function on $W_i\cap W_j$ of the line bundle associated to the divisor $\mathcal{D}_V$ and by $J_{ij}$ the Jacobian matrix of change of coordinates. Let us define the following matrices:
\begin{equation*}
 L_{ij}:=\left(\begin{array}{cc}
l_{ij}I_s & 0_{s,t-s} \\
0_{t-s,s} & I_{t-s}
\end{array}\right)
\end{equation*}
and
\begin{equation*}
 F_i:=\left(\begin{array}{cc}
f_i I_s & 0_{s,t-s} \\
0_{t-s,s} & I_{t-s}
\end{array}\right)
\end{equation*}
where $I_k$ is the identity matrix of order $k$ and $0_{k,l}$ is the $k$-by-$l$ zero matrix.

The matrix which expresses a change of basis in $\jfr(\mathcal{E})_x$ under a change of basis in $E_x$ is:
\begin{equation*}
 \left(\begin{array}{cc}
L_{ij}\otimes g_{ij} & (F_i^{-1}\otimes \mathrm{id})\cdot dg_{ij} \\
 0 & g_{ij}
\end{array}\right),
\end{equation*}
where the block at the position (1,2) is a regular matrix function, because $g_{ij}$ is constant along $\mathcal{D}_V.$ The matrix corresponding to a change of local coordinates is:
\begin{equation*}
 \left(\begin{array}{cc}
L_{ij}\cdot J_{ij}\otimes \mathrm{id} & 0 \\
0 & \mathrm{id}
\end{array}\right).
\end{equation*}
The framed sheaf of first jets $\jfr(\mathcal{E})$ of $\mathcal{E}$ fits into an exact sequence of coherent sheaves of $\mathcal{O}_{\mathcal{X}}$-modules:
\begin{equation}\label{eq:framedatiyahclass}
0\longrightarrow  \left(p_S^*(\Omega^1_S)(-\mathcal{D})\oplus p_X^*(\Omega_X^1)\right)\otimes E \longrightarrow \jfr(\mathcal{E}) \longrightarrow E\longrightarrow 0.
\end{equation}
\begin{definition}
Let $\mathcal{E}=(E,\alpha)$ be a flat family of framed sheaves parametrized by a scheme $S.$ We call \textit{framed Atiyah class} of the family $\mathcal{E}$ the class $at(\mathcal{E})$ in 
\begin{equation*}
\mathrm{Ext}^1(E,\left(p_S^*(\Omega^1_S)(-\mathcal{D})\oplus p_X^*(\Omega_X^1)\right)\otimes E) 
\end{equation*}
associated to the extension \eqref{eq:framedatiyahclass}.
\end{definition}
Let us consider the short exact sequence
\begin{equation*}
0\longrightarrow p_S^{*}(\Omega_S^1)(-\mathcal{D})\oplus p_X^*(\Omega_X^1)\stackrel{i}{\longrightarrow} \Omega^1_{\mathcal{X}}\stackrel{q}{\longrightarrow} p_S^{*}(\Omega_S^1)\vert_{\mathcal{D}}\longrightarrow 0. 
\end{equation*}
After tensoring by $E$ and applying the functor $\mathrm{Hom}(E,\cdot)$, we get the long exact sequence
\begin{eqnarray*}
&&\cdots\rightarrow\mathrm{Ext}^1(E, \left(p_S^{*}(\Omega_S^1)(-\mathcal{D})\oplus p_X^*(\Omega_X^1)\right)\otimes E)\stackrel{i_*}{\longrightarrow}\\
&&\mathrm{Ext}^1(E,\Omega^1_{\mathcal{X}}\otimes E)\stackrel{q_*}{\longrightarrow}\mathrm{Ext}^1(E, p_S^{*}(\Omega_S^1)\vert_{\mathcal{D}}\otimes E)\rightarrow\cdots.
\end{eqnarray*}
By construction, the image of $at(\mathcal{E})$ under $i_*$ is $at(E)$, which is equivalent to saying that we have the commutative diagram
\begin{equation*}
  \begin{tikzpicture}[xscale=3.2,yscale=-1.2]
    \node (A0_0) at (0, 0) {$0$};
    \node (A0_1) at (1, 0) {$\left(p_S^*(\Omega^1_S)(-\mathcal{D})\oplus p_X^*(\Omega_X^1)\right)\otimes E$};
    \node (A0_2) at (2.5, 0) {$\jfr(\mathcal{E})$};
    \node (A0_3) at (3.5, 0) {$E$};
    \node (A0_4) at (4, 0) {$0$};
    \node (A1_0) at (0, 1) {$0$};
    \node (A1_1) at (1, 1) {$\Omega_{\mathcal{X}}^1\otimes E$};
    \node (A1_2) at (2.5, 1) {$\jf(E)$};
    \node (A1_3) at (3.5, 1) {$E$};
    \node (A1_4) at (4, 1) {$0$};
    \path (A0_1) edge [->]node [auto] {$\scriptstyle{}$} (A1_1);
    \path (A0_0) edge [->]node [auto] {$\scriptstyle{}$} (A0_1);
    \path (A0_1) edge [->]node [auto] {$\scriptstyle{}$} (A0_2);
    \path (A1_0) edge [->]node [auto] {$\scriptstyle{}$} (A1_1);
    \path (A0_3) edge [thin, double distance=1.5pt]node [auto] {$\scriptstyle{}$} (A1_3);
    \path (A1_1) edge [->]node [auto] {$\scriptstyle{}$} (A1_2);
    \path (A0_3) edge [->]node [auto] {$\scriptstyle{}$} (A0_4);
    \path (A0_2) edge [->]node [auto] {$\scriptstyle{}$} (A1_2);
    \path (A1_2) edge [->]node [auto] {$\scriptstyle{}$} (A1_3);
    \path (A0_2) edge [->]node [auto] {$\scriptstyle{}$} (A0_3);
    \path (A1_3) edge [->]node [auto] {$\scriptstyle{}$} (A1_4);
  \end{tikzpicture} 
\end{equation*}
Moreover, $q_*(at(E))=q_*(i_*(at(\mathcal{E})))=0$, hence we get the commutative diagram
\begin{equation*}
  \begin{tikzpicture}[xscale=3.2,yscale=-1.2]
    \node (A0_0) at (0, 0) {$0$};
    \node (A0_1) at (0.7, 0) {$\Omega_{\mathcal{X}}^1\otimes E$};
    \node (A0_2) at (2, 0) {$\jf(E)$};
    \node (A0_3) at (3.5, 0) {$E$};
    \node (A0_4) at (4, 0) {$0$};
    \node (A1_0) at (0, 1) {$0$};
    \node (A1_1) at (0.7, 1) {$p_S^{*}(\Omega_S^1)\vert_{\mathcal{D}}\otimes E$};
    \node (A1_2) at (2, 1) {$\left(p_S^{*}(\Omega_S^1)\vert_{\mathcal{D}}\otimes E\right)\oplus E$};
    \node (A1_3) at (3.5, 1) {$E$};
    \node (A1_4) at (4, 1) {$0$};
    \path (A0_1) edge [->]node [auto] {$\scriptstyle{}$} (A1_1);
    \path (A0_0) edge [->]node [auto] {$\scriptstyle{}$} (A0_1);
    \path (A0_1) edge [->]node [auto] {$\scriptstyle{}$} (A0_2);
    \path (A1_0) edge [->]node [auto] {$\scriptstyle{}$} (A1_1);
    \path (A0_3) edge [thin, double distance=1.5pt]node [auto] {$\scriptstyle{}$} (A1_3);
    \path (A1_1) edge [->]node [auto] {$\scriptstyle{}$} (A1_2);
    \path (A0_3) edge [->]node [auto] {$\scriptstyle{}$} (A0_4);
    \path (A0_2) edge [->]node [auto] {$\scriptstyle{}$} (A1_2);
    \path (A1_2) edge [->]node [auto] {$\scriptstyle{}$} (A1_3);
    \path (A0_2) edge [->]node [auto] {$\scriptstyle{}$} (A0_3);
    \path (A1_3) edge [->]node [auto] {$\scriptstyle{}$} (A1_4);
  \end{tikzpicture} 
\end{equation*}
\begin{example}\label{ex:linebundle}
Let $F$ be a line bundle on $D.$ Let $\mathcal{L}=(L, \alpha)$ be a flat family of framed sheaves parametrized by $S$ with $L$ line bundle on $\mathcal{X}.$ In this case $V$ is the whole $\mathcal{X}$ and we choose transition functions $g_{ij}^0$ and $g_{ij}$ for $p_X^*(F_D)$ and $L$, respectively, such that
\begin{equation*}
g_{ij}\vert_{\mathcal{D}}= g_{ij}^0. 
\end{equation*}
Recall that $dg_{ij} g_{ij}^{-1}$ is a cocycle representing $at(L)$ (cf. \cite[Proposition 12]{art:atiyah1957}). By the choice of $g_{ij}^0$, we get that $d_S(g_{ij})$ vanishes along $\mathcal{D}$, where $d_S$ is the exterior differential of $S.$ Hence $dg_{ij} g_{ij}^{-1}$ can be also interpreted as a cocycle representing $at(\mathcal{L}).$ Moreover, it vanishes under the restriction of the de Rham differential $\tilde{d}:=d\vert_{p_S^*(\Omega^1_S)(-\mathcal{D})\oplus p_X^*(\Omega_X^1)}.$
\triend
\end{example}
Now we shall provide another way to describe the framed Atiyah class of a flat family of framed sheaves $\mathcal{E}=(E,\alpha)$ by using finite locally free resolutions of $E$, but in this case the costruction is \emph{local over the base}, as we will explain in the following. First, we recall a result due to B{\u{a}}nic{\u{a}}, Putinar and Schumacher that will be very useful later on. 
\begin{theorem}{\normalfont (\cite[Satz 3]{art:banicaputinarschumacher1980}).}\label{thm:banica}
Let $p\colon R\rightarrow T$ be a flat proper morphism of schemes of finite type over $k$, $T$ smooth, $E$ and $G$ coherent $\mathcal{O}_R$-modules, flat over $T.$ If the function $y\mapsto \dim \mathrm{Ext}^l(E_y,G_y)$ is constant for $l$ fixed, then the sheaf $\mathcal{E}xt_p^l(E,G)$ is locally free on $T$ and for any $y\in T$ we have
\begin{equation*}
\mathcal{E}xt_p^i(E,G)_y\otimes_{\mathcal{O}_{T,y}}\left(\sfrac{\mathcal{O}_{T,y}}{m_y}\right)\cong  \mathrm{Ext}^i(E_y,G_y)\, \mbox{ for }\, i=l-1, l.
\end{equation*}
Moreover, the same statement is true for complexes.
\end{theorem}
Let $\mathcal{E}=(E,\alpha)$ be a flat family of framed sheaves parametrized by a smooth scheme $S.$ Since the projection morphism $p_S\colon \mathcal{X}\longrightarrow S$ is smooth and projective, there exists a finite locally free resolution $E^\bullet \rightarrow E$ of $E$ (see, e.g., \cite[Proposition 2.1.10]{book:huybrechtslehn2010}).

Let us fix a point $s_0\in S.$ By the flatness property, the complex $(E^\bullet)\vert_{\{s_0\}\times D}$ is a finite resolution of locally free sheaves of $E\vert_{\{s_0\}\times D}\cong F_D.$ Let us denote by $F^\bullet$ the complex $(E^\bullet)\vert_{\{s_0\}\times D}.$ Define $\mathcal{F}^\bullet:=F^\bullet \boxtimes \mathcal{O}_S.$ The complex $\mathcal{F}^\bullet$ is $S$-flat since $(E^\bullet)\vert_{\{s_0\}\times D}$ is a complex of locally free $\mathcal{O}_D$-modules and the sheaf $\mathcal{O}_{\mathcal{D}}$ is a $S$-flat $\mathcal{O}_{\mathcal{X}}$-module. Moreover, for any $s\in S$, the complex $(\mathcal{F}^\bullet)\vert_{\{s\}\times X}$ is quasi-isomorphic to $F$, hence we get
\begin{equation*}
\mathrm{Hom}((E^\bullet)\vert_{\{s\}\times X}, (\mathcal{F}^\bullet)\vert_{\{s\}\times X})=\mathrm{Hom}(E\vert_{\{s\}\times X}, F_D)\cong\mathrm{End}(F_D).
\end{equation*}
By applying Theorem \ref{thm:banica}, we get that the natural morphism of complexes between $E^\bullet$ and $\mathcal{F}^\bullet$ on $\{s_0\}\times X$ extends to a morphism of complexes
\begin{equation*}
\alpha_\bullet \colon E^\bullet\longrightarrow\mathcal{F}^\bullet.
\end{equation*}
Let $U\subset S$ be a neighborhood of $s_0$ such that the following condition holds
\begin{equation*}
(\alpha_\bullet)\vert_{\{s\}\times D}\, \mbox{ \emph{is an isomorphism for any} } s\in U.
\end{equation*}
Let $\mathcal{X}_U=U\times X$ and $\mathcal{D}_U=U\times D.$ For any $i$, the pair $\mathcal{E}_U^i:=(E^i\vert_{\mathcal{X}_U},\alpha_i\vert_{\mathcal{X}_U}\colon E^i\vert_{\mathcal{X}_U} \rightarrow\mathcal{F}^i\vert_{\mathcal{X}_U})$ is a flat family $\mathcal{E}_U^i$ of $(D,E^i\vert_{\{s_0\}\times D})$-framed sheaves parametrized by $U.$ 

Thus we proved the following:
\begin{proposition}\label{prop:locallyfree}
Let $\mathcal{E}=(E,\alpha)$ be a flat family of framed sheaves parametrized by a smooth scheme $S$ and $E^\bullet \rightarrow E$ a finite locally free resolution of $E.$ Let $s_0$ be a point in $S.$ Then there exists a complex $\mathcal{F}^\bullet$, a morphism of complexes $\alpha_\bullet\colon E^\bullet \rightarrow \mathcal{F}^\bullet$ and an open neighborhood $U\subset S$ of $s_0$ with the following property: for any $i$ the sheaf $\mathcal{F}^i\vert_{\{s_0\}\times D}$ is a locally free $\mathcal{O}_D$-module and the pair $\mathcal{E}_U^i:=(E^i\vert_{\mathcal{X}_U}, \alpha_i\vert_{\mathcal{X}_U})$ is a flat family of $(D,\mathcal{F}^i\vert_{\{s_0\}\times D})$-framed sheaves parametrized by $U.$
\end{proposition}
If for any $i$, we consider the short exact sequence associated to $\jfr(\mathcal{E}_U^i)$
\begin{equation*}
0\longrightarrow  \left(p_U^*(\Omega^1_U)(-\mathcal{D}_U)\oplus p_X^*(\Omega_X^1)\right)\otimes E^i\vert_{\mathcal{X}_U} \longrightarrow \jfr(\mathcal{E}_U^i) \longrightarrow E^i\vert_{\mathcal{X}_U}\longrightarrow 0,
\end{equation*}
we get a class $at_{U}(\mathcal{E})$ in
\begin{eqnarray*}
&&\HE^1\left(E^\bullet\vert_{\mathcal{X}_U},\left(p_U^*(\Omega^1_U)(-\mathcal{D}_U)\oplus p_X^*(\Omega_X^1)\right)\otimes E^\bullet\vert_{\mathcal{X}_U}\right)\cong\\
&&\cong \mathrm{Ext}^1\left(E\vert_{\mathcal{X}_U},\left(p_U^*(\Omega^1_U)(-\mathcal{D}_U)\oplus p_X^*(\Omega_X^1)\right)\otimes E\vert_{\mathcal{X}_U}\right).
\end{eqnarray*}
By construction, $at_{U}(\mathcal{E})$ is independent of the resolution and is the image of $at(\mathcal{E})$ with respect to the map on Ext-groups induced by the natural morphism $i^*\colon \Omega_S^1\rightarrow \Omega_U^1$, where $i\colon U\hookrightarrow S$ is the inclusion morphism.

\subsection{Framed Newton polynomials}

Let $\mathcal{E}=(E,\alpha)$ be a flat family of framed sheaves parametrized by a smooth scheme $S.$ Let $at(\mathcal{E})^i$ denote the image in $\mathrm{Ext}^i\left(E,\tilde{\Omega}_{\mathcal{X}}^i\otimes E\right)$ of the $i$-th product
\begin{equation*}
 at(\mathcal{E})\circ \cdots \circ at(\mathcal{E})\in \mathrm{Ext}^i\left(E,(\tilde{\Omega}_{\mathcal{X}}^1)^{\otimes i}\otimes E\right)
\end{equation*}
under the morphism induced by $(\tilde{\Omega}_{\mathcal{X}}^1)^{\otimes i}\rightarrow \tilde{\Omega}_{\mathcal{X}}^i$, where $\tilde{\Omega}_{\mathcal{X}}^1:=p_S^*(\Omega^1_S)(-\mathcal{D})\oplus p_X^*(\Omega_X^1)$ and $\tilde{\Omega}_{\mathcal{X}}^i:=\Lambda^i(\tilde{\Omega}_{\mathcal{X}}^1)$ is the $i$-th exterior power of $\tilde{\Omega}_{\mathcal{X}}^1.$
\begin{definition}
The $i$-th \emph{framed Newton polynomial} of $\mathcal{E}$ is
\begin{equation*}
 \gamma^i(\mathcal{E}):=\mathrm{tr}(at(\mathcal{E})^i)\in \mathrm{H}^i(\mathcal{X},\tilde{\Omega}_{\mathcal{X}}^i).
\end{equation*}
\end{definition}
Fix a finite locally free resolution $E^\bullet\rightarrow E$ of $E$ and a point $s_0\in S.$ Let $U\subset S$ be a neighborhood of $s_0$ as in Proposition \ref{prop:locallyfree}. The restriction of $\gamma^i(\mathcal{E})$ to $\mathcal{X}_U$ coincides with the class $\gamma_{U}^i(\mathcal{E})$ defined as
\begin{equation*}
 \gamma_{U}^i(\mathcal{E}):=\mathrm{tr}(at_{U}(\mathcal{E})^i)\in \mathrm{H}^i(\mathcal{X}_U,\tilde{\Omega}_{\mathcal{X}_U}^i).
\end{equation*}

The restricted de Rham differential $\tilde{d}$ introduced in Example \ref{ex:linebundle} induces $k$-linear maps 
\begin{equation*}
 \tilde{d}\colon \mathrm{H}^i(\mathcal{X},\tilde{\Omega}_{\mathcal{X}}^{i})\longrightarrow \mathrm{H}^{i}(\mathcal{X},\tilde{\Omega}_{\mathcal{X}}^{i+1}(\mathcal{D})).
\end{equation*}
\begin{proposition}\label{prop:tildedclosed}
The $i$-th \emph{framed Newton polynomial} of $\mathcal{E}$ is $\tilde{d}$-closed. 
\end{proposition}
\proof
Let $E^\bullet\rightarrow E$ be a finite locally free resolution of $E$ and $s_0\in S.$ Let $U\subset S$ be a neighborhood of $s_0$ as in Proposition \ref{prop:locallyfree}. The cohomology class $\gamma_{U}^i(\mathcal{E})$ is $\tilde{d}\vert_U$-closed by the same arguments as in the nonframed case (cf. \cite[Section 10.1.6]{book:huybrechtslehn2010}), in particular the fact that we can reduce to the case of line bundles by using the splitting principle, and by Example \ref{ex:linebundle}. Since the restriction of $\gamma^i(\mathcal{E})$ to $\mathcal{X}_U$ is $\gamma_{U}^i(\mathcal{E})$, we get that $\gamma^i(\mathcal{E})$ is closed with respect to $\tilde{d}.$
\endproof

\subsection{The Kodaira-Spencer map for framed sheaves}

Let $\mathcal{E}=(E,\alpha)$ be a flat family of framed sheaves parametrized by a scheme $S.$ Consider the framed Atiyah class $at(\mathcal{E})$ in $\mathrm{Ext}^1\left(E,\left(p_S^*(\Omega^1_S)(-\mathcal{D})\oplus p_X^*(\Omega_X^1)\right)\otimes E\right)$ and the induced section $\mathcal{A}t(\mathcal{E})$ under the relative local-to-global map
\begin{equation*}
\mathrm{Ext}^1\left(E,\left(p_S^*(\Omega^1_S)(-\mathcal{D})\oplus p_X^*(\Omega_X^1)\right)\otimes E\right)\longrightarrow \mathrm{H}^0(S, \mathcal{E}xt_{p_S}^1(E,\left(p_S^*(\Omega^1_S)(-\mathcal{D})\oplus p_X^*(\Omega_X^1)\right)\otimes E)),
\end{equation*}
coming from the relative local-to-global spectral sequence 
\begin{equation*}
\mathrm{H}^i(S,\mathcal{E}xt_{p_S}^j(E,\left(p_S^*(\Omega^1_S)(-\mathcal{D})\oplus p_X^*(\Omega_X^1)\right)\otimes E))\Rightarrow \mathrm{Ext}^{i+j}(E,\left(p_S^*(\Omega^1_S)(-\mathcal{D})\oplus p_X^*(\Omega_X^1)\right)\otimes E). 
\end{equation*}
By considering the $S$-part $\mathcal{A}t_S(\mathcal{E})$ of $\mathcal{A}t(\mathcal{E})$ in $\mathrm{H}^0(S,\mathcal{E}xt_{p_S}^1(E,p_S^*(\Omega^1_S)(-\mathcal{D})\otimes E))$, we define the framed version of the Kodaira-Spencer map.
\begin{definition}
The \textit{framed Kodaira-Spencer map} associated to the family $\mathcal{E}$ is the composition
\begin{eqnarray*}
 KS_{fr}\colon (\Omega_S^1)^\vee&\stackrel{\mathrm{id}\otimes \mathcal{A}t_S(\mathcal{E})}{\longrightarrow}& (\Omega_S^1)^\vee\otimes \mathcal{E}xt_{p_S}^1(E,p_S^*(\Omega^1_S)(-\mathcal{D})\otimes E)\rightarrow\\
&\longrightarrow& \mathcal{E}xt_{p_S}^1(E,p_S^*((\Omega_S^1)^\vee\otimes\Omega^1_S)\otimes p_X^*(\mathcal{O}_X(-D))\otimes E)\rightarrow\\
&\longrightarrow& \mathcal{E}xt_{p_S}^1(E,p_X^*(\mathcal{O}_X(-D))\otimes E).
\end{eqnarray*}
\end{definition}

\subsection{Closed two-forms via the framed Atiyah class}

Let $S$ be a smooth affine scheme and $\mathcal{E}=(E,\alpha)$ a flat family of framed sheaves parametrized by $S.$ Let $\gamma^{0,2}$ denote the component of $\gamma^2(\mathcal{E})$ in $\mathrm{H}^0(S, \Omega^2_S)\otimes \mathrm{H}^{2}(X,\mathcal{O}_X(-2\,D)).$
\begin{definition}
Let $\tau_S$ be the homomorphism given by
\begin{equation*}
\tau_S\colon \mathrm{H}^0(X,\omega_X(2\,D))\cong \mathrm{H}^2(X,\mathcal{O}_X(-2\,D))^\vee\stackrel{\cdot \,\gamma^{0,2}}{\longrightarrow} \mathrm{H}^0(S, \Omega^2_S),
\end{equation*}
 where $\cong$ denotes Serre's duality.
\end{definition}
\begin{proposition}\label{prop:closedness}
For any $\omega\in \mathrm{H}^0(X,\omega_X(2\,D))$, the associated two-form $\tau_S(\omega)$ on $S$ is closed.
\end{proposition}
\proof
We can write
\begin{equation*}
\gamma^{0,2}=\sum_l \mu_l\otimes \nu_l,
\end{equation*}
for elements $\mu_l\in \mathrm{H}^0(S, \Omega^2_S)$ and $\nu_l\in \mathrm{H}^2(X,\mathcal{O}_X(-2\,D)).$ Since $\tilde{d}(\gamma^2(\mathcal{E}))=0$ (cf. Proposition \ref{prop:tildedclosed}), the component of $\tilde{d}(\gamma^{0,2})$ in $\mathrm{H}^0(S, \Omega^3_S)\otimes \mathrm{H}^{2}(X,\mathcal{O}_X(-2\,D))$ is zero, which means
\begin{equation*}
\sum_l d_{S}(\mu_l)\otimes \nu_l=0.
\end{equation*}
Therefore
\begin{equation*}
d_S(\tau_S(\omega))=d_S\left(\sum_l \mu_l \cdot \omega(\nu_l)\right)=\sum_l d_{S}(\mu_l)\cdot\omega(\nu_l)=0.
\end{equation*}
\endproof
Fix $\omega\in \mathrm{H}^0(X,\omega_X(2\,D)).$  Since $S$ is smooth, it follows from the definitions of the framed Kodaira-Spencer map and the framed Newton polynomial that the two-form $\tau_S(\omega)$ at a point $s_0\in S$ coincides with the following composition of maps:
\begin{eqnarray*}
&&T_{s_0} S\times T_{s_0} S\stackrel{KS \times KS}{\longrightarrow}\mathrm{Ext}^1(E\vert_{\{s_0\}\times X},E\vert_{\{s_0\}\times X}(-D))\times \mathrm{Ext}^1(E\vert_{\{s_0\}\times X},E\vert_{\{s_0\}\times X}(-D))\\
&& \stackrel{\circ}{\longrightarrow} \mathrm{Ext}^2(E\vert_{\{s_0\}\times X},E\vert_{\{s_0\}\times X}(-2\,D))\stackrel{tr}{\longrightarrow} \mathrm{H}^2(X, \mathcal{O}_X(-2\,D))\stackrel{\cdot\, \omega}{\longrightarrow} \mathrm{H}^2(X,\omega_X)\cong k.
\end{eqnarray*}

\section{The tangent bundle of moduli spaces of framed sheaves}\label{sec:framedkodaira}

Let $\mathcal{M}^s(X;P)$ be the moduli space of Gieseker-stable torsion-free sheaves on $X$ with Hilbert polynomial $P.$ The open subset $\mathcal{M}_0(X;P)\subset \mathcal{M}^s(X;P)$ of points $[E]$ such that $\mathrm{Ext}_0^2(E,E)$ vanishes is smooth according to \cite[Theorem 4.5.4]{book:huybrechtslehn2010}. Suppose there exists a universal family $\tilde{E}$ on $\mathcal{M}_0(X;P)\times X.$ By using the Atiyah class of $\tilde{E}$, one can define the Kodaira-Spencer map for $\mathcal{M}_0(X;P)$:
\begin{equation*}
KS\colon T\mathcal{M}_0(X;P) \longrightarrow \mathcal{E}xt_p^1(\tilde{E},\tilde{E}),
\end{equation*}
where $p\colon \mathcal{M}_0(X;P)\times X\rightarrow \mathcal{M}_0(X;P)$ is the projection. Moreover, it is possible to prove that $KS$ is an isomorphism (this result holds also when a universal family for $\mathcal{M}_0(X;P)$ does not exist, cf. \cite[Theorem 10.2.1]{book:huybrechtslehn2010}). In this section we shall prove the framed analogue of this result for the moduli spaces of stable framed sheaves on $X.$

Let $\delta\in \mathbb{Q}[n]$ be a stability polynomial and $P$ a numerical polynomial of degree two. Let $\mathcal{M}^*_\delta(X;F_D,P)$ be the moduli space of framed sheaves on $X$ with Hilbert polynomial $P$ that are stable with respect to $\mathcal{O}_X(1)$ and $\delta.$ This is an open subset of the fine moduli space $\mathcal{M}_\delta^{s}(X;F_D,P)$ of stable $F_D$-framed modules with Hilbert polynomial $P.$ 

Let us denote by $\mathcal{M}^*_\delta(X;F_D,P)^{sm}$ the smooth locus of $\mathcal{M}^*_\delta(X;F_D,P)$ and by $\tilde{\mathcal{E}}=(\tilde{E},\tilde{\alpha})$ a \emph{universal framed sheaf} over $\mathcal{M}^*_\delta(X;F_D,P)^{sm}.$ 
\begin{theorem}\label{thm:tanbundle}
The framed Kodaira-Spencer map defined by $\tilde{\mathcal{E}}$ induces a canonical isomorphism
\begin{equation*}
KS_{fr}\colon T\mathcal{M}^*_\delta(X;F_D,P)^{sm}\stackrel{\sim}{\longrightarrow} \mathcal{E}xt_p^1(\tilde{E}, \tilde{E}\otimes p_X^*(\mathcal{O}_X(-D))),
\end{equation*}
where $p$ is the projection from $\mathcal{M}^*_\delta(X;F_D,P)^{sm}\times X$ to $\mathcal{M}^*_\delta(X;F_D,P)^{sm}.$
\end{theorem}
\proof
First note that $\mathcal{M}^*_\delta(X;F_D,P)^{sm}$ is a reduced separated scheme of finite type over $k.$ Hence it suffices to prove that the framed Kodaira-Spencer map is an isomorphism on the fibres over closed points. 

Let $[(E,\alpha)]$ be a closed point. We want to show that the Kodaira-Spencer map $KS_{fr}([(E,\alpha)])$ on the fibre over $[(E,\alpha)]$ coincides with the isomorphism
\begin{equation*}
T_{[(E,\alpha)]} \mathcal{M}^*_\delta(X;F_D,P)^{sm}\stackrel{\sim}{\longrightarrow} \mathrm{Ext}^1(E,E(-D)),
\end{equation*} 
coming from deformation theory (see proof of \cite[Theorem 4.1]{art:huybrechtslehn1995-II}).

Let $w$ be an element in $\mathrm{Ext}^1(E,E(-D)).$ Consider the long exact sequence
\begin{equation*}
 \cdots \rightarrow \mathrm{Ext}^1(E,E(-D))\stackrel{j_*}{\longrightarrow} \mathrm{Ext}^1(E,E) \stackrel{\alpha_*}{\longrightarrow} \mathrm{Ext}^1(E,F_D)\rightarrow \cdots
\end{equation*}
obtained by applying the functor $\mathrm{Hom}(E,\cdot)$ to the exact sequence
\begin{equation*}
0\longrightarrow E(-D)\stackrel{j}{\longrightarrow} E\stackrel{\alpha}{\longrightarrow} F_D\longrightarrow 0.
\end{equation*}
Let $v=j_*(w)\in \mathrm{Ext}^1(E,E).$ We get a commutative diagram
\begin{equation*}
  \begin{tikzpicture}[xscale=1.7,yscale=-1.2]
    \node (A0_0) at (0, 0) {$0$};
    \node (A0_1) at (1, 0) {$E(-D)$};
    \node (A0_2) at (2, 0) {$\tilde{G}$};
    \node (A0_3) at (3, 0) {$E$};
    \node (A0_4) at (4, 0) {$0$};
    \node (A1_0) at (0, 1) {$0$};
    \node (A1_1) at (1, 1) {$E$};
    \node (A1_2) at (2, 1) {$G$};
    \node (A1_3) at (3, 1) {$E$};
    \node (A1_4) at (4, 1) {$0$};
    \path (A0_1) edge [->]node [auto] {$\scriptstyle{j}$} (A1_1);
    \path (A0_0) edge [->]node [auto] {$\scriptstyle{}$} (A0_1);
    \path (A0_1) edge [->]node [auto] {$\scriptstyle{\tilde{i}}$} (A0_2);
    \path (A1_0) edge [->]node [auto] {$\scriptstyle{}$} (A1_1);
    \path (A0_3) edge [thin, double distance=1.5pt]node [auto] {$\scriptstyle{}$} (A1_3);
    \path (A1_1) edge [->]node [auto] {$\scriptstyle{i}$} (A1_2);
    \path (A0_3) edge [->]node [auto] {$\scriptstyle{}$} (A0_4);
    \path (A0_2) edge [->]node [auto] {$\scriptstyle{}$} (A1_2);
    \path (A1_2) edge [->]node [auto] {$\scriptstyle{\pi}$} (A1_3);
    \path (A0_2) edge [->]node [auto] {$\scriptstyle{\tilde{\pi}}$} (A0_3);
    \path (A1_3) edge [->]node [auto] {$\scriptstyle{}$} (A1_4);
  \end{tikzpicture} 
\end{equation*}
where the first row is a representative for $w$ and the second one a representative for $v.$

Let $S=\mathrm{Spec}(k[\varepsilon])$ be the spectrum of the ring of dual numbers, where $\varepsilon^2=0.$ We can think of $G$ as a $S$-flat family by letting $\varepsilon$ act on $G$ as the morphism $i\circ \pi.$ Since $\varepsilon \tilde{G}=E(-D)$ and $\varepsilon G=E$, by applying the snake lemma to the previous diagram we get
\begin{equation*}
  \begin{tikzpicture}[xscale=1.7,yscale=-1.2]
    \node (A0_0) at (0, 0) {$0$};
    \node (A0_1) at (1, 0) {$E(-D)$};
    \node (A0_2) at (2, 0) {$\tilde{G}$};
    \node (A0_3) at (3, 0) {$E$};
    \node (A0_4) at (4, 0) {$0$};
    \node (A1_0) at (0, 1) {$0$};
    \node (A1_1) at (1, 1) {$E$};
    \node (A1_2) at (2, 1) {$G$};
    \node (A1_3) at (3, 1) {$E$};
    \node (A1_4) at (4, 1) {$0$};
    \node (A2_0) at (0, 2) {$0$};
    \node (A2_1) at (1, 2) {$\varepsilon F_D$};
    \node (A2_2) at (2, 2) {$\varepsilon F_D$};
    \node (A2_3) at (3, 2) {$0$};
    \node (A3_1) at (1, 3) {$0$};
    \node (A3_2) at (2, 3) {$0$};
    \path (A2_1) edge [->]node [auto] {$\scriptstyle{}$} (A3_1);
    \path (A0_1) edge [->]node [auto] {$\scriptstyle{}$} (A1_1);
    \path (A2_1) edge [->]node [auto] {$\scriptstyle{}$} (A2_2);
    \path (A0_0) edge [->]node [auto] {$\scriptstyle{}$} (A0_1);
    \path (A0_1) edge [->]node [auto] {$\scriptstyle{\tilde{i}}$} (A0_2);
    \path (A1_0) edge [->]node [auto] {$\scriptstyle{}$} (A1_1);
    \path (A0_3) edge [->]node [auto] {$\scriptstyle{\mathrm{id}_E}$} (A1_3);
    \path (A1_1) edge [->]node [auto] {$\scriptstyle{i}$} (A1_2);
  \path (A1_1) edge [->]node [auto] {$\scriptstyle{}$} (A2_1);
  \path (A1_2) edge [->]node [auto] {$\scriptstyle{\beta}$} (A2_2);
    \path (A2_2) edge [->]node [auto] {$\scriptstyle{}$} (A2_3);
    \path (A0_3) edge [->]node [auto] {$\scriptstyle{}$} (A0_4);
    \path (A2_2) edge [->]node [auto] {$\scriptstyle{}$} (A3_2);
    \path (A0_2) edge [->]node [auto] {$\scriptstyle{}$} (A1_2);
    \path (A2_0) edge [->]node [auto] {$\scriptstyle{}$} (A2_1);
    \path (A1_2) edge [->]node [auto] {$\scriptstyle{\pi}$} (A1_3);
    \path (A1_3) edge [->]node [auto] {$\scriptstyle{}$} (A2_3);
    \path (A0_2) edge [->]node [auto] {$\scriptstyle{\tilde{\pi}}$} (A0_3);
    \path (A1_3) edge [->]node [auto] {$\scriptstyle{}$} (A1_4);
  \end{tikzpicture}
\end{equation*}
Moreover $\alpha_*(v)=0$, hence we have the commutative diagram
\begin{equation*}
  \begin{tikzpicture}[xscale=1.7,yscale=-1.2]
    \node (A0_0) at (0, 0) {$0$};
    \node (A0_1) at (1, 0) {$E$};
    \node (A0_2) at (2, 0) {$G$};
    \node (A0_3) at (3, 0) {$E$};
    \node (A0_4) at (4, 0) {$0$};
    \node (A1_0) at (0, 1) {$0$};
    \node (A1_1) at (1, 1) {$\varepsilon F_D$};
    \node (A1_2) at (2, 1) {$E\oplus\varepsilon F_D$};
    \node (A1_3) at (3, 1) {$E$};
    \node (A1_4) at (4, 1) {$0$};
    \path (A0_1) edge [->]node [auto] {$\scriptstyle{}$} (A1_1);
    \path (A0_0) edge [->]node [auto] {$\scriptstyle{}$} (A0_1);
    \path (A0_1) edge [->]node [auto] {$\scriptstyle{i}$} (A0_2);
    \path (A1_0) edge [->]node [auto] {$\scriptstyle{}$} (A1_1);
    \path (A0_3) edge [thin, double distance=1.5pt]node [auto] {$\scriptstyle{}$} (A1_3);
    \path (A1_1) edge [->]node [auto] {$\scriptstyle{}$} (A1_2);
    \path (A0_3) edge [->]node [auto] {$\scriptstyle{}$} (A0_4);
    \path (A0_2) edge [->]node [auto] {$\scriptstyle{}$} (A1_2);
    \path (A1_2) edge [->]node [auto] {$\scriptstyle{}$} (A1_3);
    \path (A0_2) edge [->]node [auto] {$\scriptstyle{\pi}$} (A0_3);
    \path (A1_3) edge [->]node [auto] {$\scriptstyle{}$} (A1_4);
  \end{tikzpicture}
\end{equation*}
Thus we obtain a framing $\gamma\colon G\rightarrow F_D\oplus \varepsilon F_D$ induced by $\alpha$ and $\beta.$ Moreover $\gamma\vert_{\mathcal{D}}$ is an isomorphism. We denote by $\mathcal{G}$ the corresponding $S$-flat family of $(D,F_D)$-framed sheaves on $X.$

In the nonframed case, one define a relative Atiyah class for families of coherent sheaves parametrized by a scheme $S$ and takes its $S$-part (see \cite[Section 10.1.8]{book:huybrechtslehn2010}). As it is explained in \cite[Example 10.1.9]{book:huybrechtslehn2010}, since $S$ is affine, the relative $S$-part $\mathcal{A}t_S(G)$ of $G$ is an element of 
\begin{equation*}
\mathrm{Ext}^1_{\mathcal{X}}(G,p_S^* \Omega^1_S\otimes G)\cong \mathrm{Ext}^1_\mathcal{X}(G,E).
\end{equation*}
Consider the short exact sequence of coherent sheaves over $\mathrm{Spec}\left(\sfrac{k[\varepsilon_1,\varepsilon_2]}{(\varepsilon_1,\varepsilon_2)^2}\right)\times X$
\begin{equation}\label{eq:s-part}
0\longrightarrow E\stackrel{i'}{\longrightarrow} G'\stackrel{\pi'}{\longrightarrow} G\longrightarrow 0,
\end{equation}
where $\varepsilon_1$ and $\varepsilon_2$ act trivially on $E$ and by $i\circ \pi$ on $G$, and $G'\cong \sfrac{k[\varepsilon_1]\otimes_{k} G}{\varepsilon_1\varepsilon_2 G}\cong G\oplus E$, with actions
\begin{equation*}
\varepsilon_1=
\left(
\begin{array}{cc}
0 & \pi\\ 
0 & 0
\end{array}
\right) \,\mbox{ and }
\varepsilon_2=
\left(
\begin{array}{cc}
i\pi & 0\\ 
0 & 0
\end{array}
\right).
\end{equation*}
By definition of Atiyah class, $\mathcal{A}t_S(G)$ is precisely the extension class of the short exact sequence \eqref{eq:s-part}, considered as a sequence of $k[\varepsilon_1]\otimes\mathcal{O}_X$-modules. 

The morphism $\pi$ induces a pull-back morphism $\pi^*\colon \mathrm{Ext}_X^1(E,E) \rightarrow\mathrm{Ext}^1_\mathcal{X}(G,E)$, which is an isomorphism. As it is proved in \cite[Example 10.1.9]{book:huybrechtslehn2010}, $\pi^*(v)=\mathcal{A}t_S(G)$, indeed we have the commutative diagram
\begin{equation*}
  \begin{tikzpicture}[xscale=1.7,yscale=-1.2]
    \node (A0_0) at (0, 0) {$0$};
    \node (A0_1) at (1, 0) {$E$};
    \node (A0_2) at (2, 0) {$G'$};
    \node (A0_3) at (3, 0) {$G$};
    \node (A0_4) at (4, 0) {$0$};
    \node (A1_0) at (0, 1) {$0$};
    \node (A1_1) at (1, 1) {$E$};
    \node (A1_2) at (2, 1) {$G$};
    \node (A1_3) at (3, 1) {$E$};
    \node (A1_4) at (4, 1) {$0$};
    \path (A0_0) edge [->]node [auto] {$\scriptstyle{}$} (A0_1);
    \path (A0_1) edge [thin, double distance=1.5pt]node [auto] {$\scriptstyle{}$} (A1_1);
    \path (A1_0) edge [->]node [auto] {$\scriptstyle{}$} (A1_1);
    \path (A0_3) edge [->]node [auto] {$\scriptstyle{\pi}$} (A1_3);
    \path (A0_1) edge [->]node [auto] {$\scriptstyle{i'}$} (A0_2);
    \path (A0_2) edge [->]node [auto] {$\scriptstyle{t'}$} (A1_2);
    \path (A1_1) edge [->]node [auto] {$\scriptstyle{i}$} (A1_2);
    \path (A1_2) edge [->]node [auto] {$\scriptstyle{\pi}$} (A1_3);
    \path (A0_2) edge [->]node [auto] {$\scriptstyle{\pi'}$} (A0_3);
    \path (A0_3) edge [->]node [auto] {$\scriptstyle{}$} (A0_4);    
    \path (A1_3) edge [->]node [auto] {$\scriptstyle{}$} (A1_4);
  \end{tikzpicture}
\end{equation*}
Thus $G'$ is the sheaf of first jets of $G$ relative to the quotient $\Omega_\mathcal{X}^1\rightarrow p_S^*(\Omega_S^1)\rightarrow 0.$ By following Maakestad's construction of Atiyah classes of coherent sheaves relative to quotients of $\Omega_\mathcal{X}^1$ (cf. \cite[Section 3]{art:maakestad2010-I}) and by readapting to this particular case our construction of the framed sheaf of first jets given in Section \ref{sec:framedatiyah}, we can define a \textit{framed sheaf of first jets} $\tilde{G}'$ of the framed sheaf $\mathcal{G}$ relative to $p_S^*(\Omega_S^1).$ Thus we get a commutative diagram
\begin{equation*}
  \begin{tikzpicture}[xscale=1.7,yscale=-1.2]
    \node (A0_0) at (0, 0) {$0$};
    \node (A0_1) at (1, 0) {$E(-D)$};
    \node (A0_2) at (2, 0) {$\tilde{G}'$};
    \node (A0_3) at (3, 0) {$G$};
    \node (A0_4) at (4, 0) {$0$};
    \node (A1_0) at (0, 1) {$0$};
    \node (A1_1) at (1, 1) {$E$};
    \node (A1_2) at (2, 1) {$G'$};
    \node (A1_3) at (3, 1) {$G$};
    \node (A1_4) at (4, 1) {$0$};
    \path (A0_1) edge [->]node [auto] {$\scriptstyle{}$} (A1_1);
    \path (A0_0) edge [->]node [auto] {$\scriptstyle{}$} (A0_1);
    \path (A0_1) edge [->]node [auto] {$\scriptstyle{\tilde{i}'}$} (A0_2);
    \path (A1_0) edge [->]node [auto] {$\scriptstyle{}$} (A1_1);
    \path (A0_3) edge [thin, double distance=1.5pt]node [auto] {$\scriptstyle{}$} (A1_3);
    \path (A1_1) edge [->]node [auto] {$\scriptstyle{i'}$} (A1_2);
    \path (A0_3) edge [->]node [auto] {$\scriptstyle{}$} (A0_4);
    \path (A0_2) edge [->]node [auto] {$\scriptstyle{}$} (A1_2);
    \path (A1_2) edge [->]node [auto] {$\scriptstyle{\pi'}$} (A1_3);
    \path (A0_2) edge [->]node [auto] {$\scriptstyle{\tilde{\pi}'}$} (A0_3);
    \path (A1_3) edge [->]node [auto] {$\scriptstyle{}$} (A1_4);
  \end{tikzpicture}
\end{equation*}
The first row is a representative for the $S$-part $\mathcal{A}t_S(\mathcal{G})$ of $\mathcal{G}$ in
\begin{equation*}
\mathrm{Ext}^1_{\mathcal{X}}(G,p_S^* \Omega^1_S(-\mathcal{D})\otimes G)\cong \mathrm{Ext}^1_\mathcal{X}(G,E(-D)).
\end{equation*}
Consider the diagram
\begin{equation*}
\begin{tikzpicture}[xscale=1.5,yscale=-1.2]
    \node (A0_0) at (0, 0) {$0$};
\node (A1_1) at (1, 1) {$0$};
\node (A0_2) at (2, 0) {$E(-D)$};
\node (A1_3) at (3, 1) {$E(-D)$};
    \node (A0_4) at (4, 0) {$\tilde{G}'$};
\node (A1_5) at (5, 1) {$\tilde{G}$};
    \node (A0_6) at (6, 0) {$G$};
\node (A1_7) at (7, 1) {$E$};
    \node (A0_8) at (8, 0) {$0$};
\node (A1_9) at (9, 1) {$0$};

\node (A2_0) at (0, 2) {$0$};
\node (A3_1) at (1, 3) {$0$};
\node (A2_2) at (2, 2) {$E$};
\node (A3_3) at (3, 3) {$E$};
\node (A2_4) at (4, 2) {$G'$};
\node (A3_5) at (5, 3) {$G$};
\node (A2_6) at (6, 2) {$G$};
\node (A3_7) at (7, 3) {$E$};
\node (A2_8) at (8, 2) {$0$};
\node (A3_9) at (9, 3) {$0$};

\node (A4_2) at (2, 4) {$F_D$};
\node (A5_3) at (3, 5) {$F_D$};
\node (A4_4) at (4, 4) {$F_D$};
\node (A5_5) at (5, 5) {$F_D$};

\path (A3_5)  edge [->>]node [auto] {$\scriptstyle{}$} (A5_5);

\path (A4_2) edge [thin, double distance=1.5pt]node [auto] {$\scriptstyle{}$} (A5_3);
\path (A5_3) edge [thin, double distance=1.5pt]node [auto] {$\scriptstyle{}$} (A5_5);
\path (A4_4) edge [thin, double distance=1.5pt]node [auto] {$\scriptstyle{}$} (A5_5);

\path (A0_0) edge [->]node [auto] {$\scriptstyle{}$} (A0_2);
\path (A0_2) edge [->]node [auto] {$\scriptstyle{\tilde{i}'}$} (A0_4);
\path (A0_4) edge [->]node [auto] {$\scriptstyle{\tilde{\pi}'}$} (A0_6);
\path (A0_6) edge [->]node [auto] {$\scriptstyle{}$} (A0_8);

\path (A1_7) edge [->]node [auto] {$\scriptstyle{}$} (A1_9);

\path (A0_2) edge [thin, double distance=1.5pt]node [auto] {$\scriptstyle{}$} (A1_3);
\path (A0_6) edge [->]node [auto] {$\scriptstyle{\pi}$} (A1_7);

\path (A2_0) edge [->]node [auto] {$\scriptstyle{}$} (A2_2);

\path (A3_5) edge [->]node [auto] {$\scriptstyle{\pi}$} (A3_7);
\path (A3_7) edge [->]node [auto] {$\scriptstyle{}$} (A3_9);

\path (A2_2) edge [thin, double distance=1.5pt]node [auto] {$\scriptstyle{}$} (A3_3);
\path (A2_4) edge [->]node [auto] {$\scriptstyle{t'}$} (A3_5);
\path (A2_6) edge [->]node [auto] {$\scriptstyle{\pi}$} (A3_7);


\path (A0_2) edge [right hook->]node [auto] {$\scriptstyle{}$} (A2_2);
\path (A1_1) edge [-,draw=white,line width=4pt]node [auto] {$\scriptstyle{}$} (A1_3);
\path (A1_1) edge [->]node [auto] {$\scriptstyle{}$} (A1_3);

\path (A0_4) edge [right hook->]node [auto] {$\scriptstyle{}$} (A2_4);
\path (A1_3) edge [-,draw=white,line width=4pt]node [auto] {$\scriptstyle{}$} (A1_5);
\path (A1_3) edge [->]node [auto] {$\scriptstyle{\tilde{i}}$} (A1_5);

\path (A0_6) edge [thin, double distance=1.5pt]node [auto] {$\scriptstyle{}$} (A2_6);
\path (A1_5) edge [-,draw=white,line width=4pt]node [auto] {$\scriptstyle{}$} (A1_7);
\path (A1_5) edge [->,near start]node [auto] {$\scriptstyle{\tilde{\pi}}$} (A1_7);

\path (A2_2)  edge [->>]node [auto] {$\scriptstyle{}$} (A4_2);
\path (A3_1) edge [-,draw=white,line width=4pt]node [auto] {$\scriptstyle{}$} (A3_3);
\path (A3_1) edge [->]node [auto] {$\scriptstyle{}$} (A3_3);

\path (A2_4)  edge [->>]node [auto] {$\scriptstyle{}$} (A4_4);
\path (A3_3) edge [-,draw=white,line width=4pt]node [auto] {$\scriptstyle{}$} (A3_5);
\path (A3_3) edge [->,near start]node [auto] {$\scriptstyle{i}$} (A3_5);

\path (A2_2) edge [->,near start]node [auto] {$\scriptstyle{i'}$} (A2_4);
\path (A1_3) edge [-,draw=white,line width=4pt]node [auto] {$\scriptstyle{}$} (A3_3);
\path (A1_3) edge [right hook->]node [auto] {$\scriptstyle{}$} (A3_3);

\path (A2_4) edge [->,near start]node [auto] {$\scriptstyle{\pi'}$} (A2_6);
\path (A1_5) edge [-,draw=white,line width=4pt]node [auto] {$\scriptstyle{}$} (A3_5);
\path (A1_5) edge [right hook->]node [auto] {$\scriptstyle{}$} (A3_5);

\path (A2_6) edge [->]node [auto] {$\scriptstyle{}$} (A2_8);
\path (A1_7) edge [-,draw=white,line width=4.5pt]node [auto] {$\scriptstyle{}$} (A3_7);
\path (A1_7) edge [thin, double distance=1.5pt]node [auto] {$\scriptstyle{}$} (A3_7);

\path (A4_2) edge [thin, double distance=1.5pt]node [auto] {$\scriptstyle{}$} (A4_4);
\path (A3_3) edge [-,draw=white,line width=4.5pt]node [auto] {$\scriptstyle{}$} (A5_3);
\path (A3_3)  edge [->>]node [auto] {$\scriptstyle{}$} (A5_3);

\end{tikzpicture}
\end{equation*}
By diagram chasing, one can define a morphism $\tilde{G}'\rightarrow \tilde{G}$ such that the corresponding diagram commutes. Thus the image of $w$ via the map $\mathrm{Ext}_X^1(E,E(-D)) \rightarrow\mathrm{Ext}^1_\mathcal{X}(G,E(-D))$ is exactly $\mathcal{A}t_S(\mathcal{G}).$ This completes the proof.
\endproof

\section{Closed two-forms on moduli spaces of framed sheaves}\label{sec:framedtwoforms}

In this section we show how to construct closed two-forms on the moduli space $\mathcal{M}^*_\delta(X;F_D,P)^{sm}$ by using global sections of the line bundle $\omega_X(2\,D).$ Moreover, we establish a nondegeneracy criterion for these two-forms.

Let $[(E,\alpha)]$ be a point in $\mathcal{M}^*_\delta(X;F_D,P)^{sm}.$ By Theorem \ref{thm:tanbundle}, the vector space $\mathrm{Ext^1}(E,E(-D))$ is naturally identified with the tangent space $T_{[(E,\alpha)]}\mathcal{M}^*_\delta(X;F_D,P).$ 

For any $\omega\in \mathrm{H}^0(X,\omega_X(2\,D))$, we can define a skew-symmetric bilinear form
\begin{eqnarray*}
&&\mathrm{Ext}^1(E,E(-D))\times \mathrm{Ext}^1(E,E(-D))\stackrel{\circ}{\longrightarrow} \mathrm{Ext}^2(E,E(-2\,D))\\
&&\stackrel{tr}{\longrightarrow} \mathrm{H}^2(X, \mathcal{O}_X(-2\,D))\stackrel{\cdot\, \omega}{\longrightarrow} \mathrm{H}^2(X,\omega_X)\cong k.
\end{eqnarray*}
By varying the point $[(E,\alpha)]$, these forms fit into a exterior two-form $\tau(\omega)$ on $\mathcal{M}^*_\delta(X;F_D,P)^{sm}.$
\begin{theorem}\label{thm:closed}
For any $\omega\in \mathrm{H}^0(X,\omega_X(2\,D))$, the two-form $\tau(\omega)$ is closed on $\mathcal{M}^*_\delta(X;F_D,P)^{sm}.$
\end{theorem}
\proof
It suffices to prove that given a smooth affine variety $S$, for any $S$-flat family $\mathcal{E}=(E,\alpha)$ of framed sheaves on $X$ defining a classifying morphism
\begin{eqnarray*}
 \psi\colon S&\longrightarrow & \mathcal{M}^*_\delta(X;F_D,P)^{sm},\\
s&\longmapsto& [\mathcal{E}\vert_{\{s\}\times X}],
\end{eqnarray*}
the pullback $\psi^*(\tau(\omega))\in \mathrm{H}^0(S,\Omega^2_S)$ is closed. Since $\psi^{*}(\tau(\omega))=\tau_S(\omega)$ by construction, this follows from Proposition \ref{prop:closedness}.
\endproof
Thus we have constructed closed two-forms $\tau(\omega)$ on the moduli space $\mathcal{M}^*_\delta(X;F_D,P)^{sm}$ depending on a choice of $\omega\in \mathrm{H}^0(X,\omega_X(2\,D)).$ In general, these forms may be degenerate. 

Now we give a criterion to check when the two-forms are non-degenerate.
\begin{proposition}\label{prop:non-dege}
Let $\omega\in \mathrm{H}^0(X,\omega_X(2\,D))$ and $[(E,\alpha)]$ a point in $\mathcal{M}^*_\delta(X;F_D,P)^{sm}.$ The closed two-form $\tau(\omega)_{[(E,\alpha)]}$ is non-degenerate at the point $[(E,\alpha)]$ if and only if the multiplication by $\omega$ induces an isomorphism
\begin{equation*}
 \omega_*\colon \mathrm{Ext}^1(E,E(-D))\longrightarrow \mathrm{Ext}^1(E,E\otimes \omega_X(D)).
\end{equation*}
\end{proposition}
\proof
The proof is similar to that of \cite[Proposition 10.4.1]{book:huybrechtslehn2010}.
\endproof
Obviously, if the line bundle $\omega_X(2\, D)$ is trivial, for any point $[(E,\alpha)]$ in $\mathcal{M}^*_\delta(X;F_D,P)^{sm}$ the pairing
\begin{equation*}
 \tau(1)\colon \mathrm{Ext}^1(E,E(-D))\times \mathrm{Ext}^1(E,E(-D))\longrightarrow k
\end{equation*}
is a non-degenerate alternating form, where $1\in \mathrm{H}^0(X,\omega_X(2\,D))\cong k.$

\section{Examples}\label{sec:example}

In this section we provide explicit examples of holomorphic symplectic structures on moduli spaces of $(D,F_D)$-framed sheaves with fixed Hilbert polynomial on some birationally ruled surfaces. In particular, if one compares our costruction of symplectic structures on the moduli spaces of framed sheaves on $\mathbb{CP}^1\times \mathbb{CP}^1$, $C\times \mathbb{CP}^1$ with $C$ elliptic curve, and on the second Hirzebruch surface $\mathbb{F}_2$ to that of Bottacin (cf. \cite{art:bottacin2000}), one can see that they are equivalent on the locally free part of the moduli space. On the other hand, we construct new examples of symplectic structures on the moduli spaces of $(D, F_D)$-framed sheaves over $\mathbb{CP}^2$ and the blowup of $\mathbb{CP}^2$ at a point, not covered by Bottacin's result.

These constructions are based upon the results of Theorem \ref{thm:bruzzomarkushevich} and Corollary \ref{cor:smoothness}. Useful references on birationally ruled surfaces are \cite[Chapter V]{book:hartshorne1977} and \cite{book:beauville1996}.

\subsection*{The complex projective plane}

Let us denote by $l_\infty$ a line in $\mathbb{CP}^2.$ Fix a positive integer $d$ and a smooth connected curve $D$ in the complete linear system $\vert d\, l_\infty\vert.$ The genus of $D$ is $\frac{(d-1)(d-2)}{2}.$ 

Let $c, n\in \mathbb{Z}$ and $F_D$ a Gieseker-semistable locally free $\mathcal{O}_D$-module of rank $r$ and degree $cd$ on $D.$ Let us denote by $\mathcal{M}^*(\mathbb{CP}^2;F_D,r,c\,l_\infty,n)$ the moduli space of  $(D,F_D)$-framed sheaves on $\mathbb{CP}^2$ with rank $r$, first Chern class $c\,l_\infty$ and second Chern class $n.$ Since $K_{\mathbb{CP}^2}=-3\, l_\infty$, we get $K_{\mathbb{CP}^2}+2\, D=(-3+2d) \,l_\infty.$ Hence $\mathrm{H}^0(\mathbb{CP}^2,\omega_{\mathbb{CP}^2}(2\,D))$ is nonzero only for $2d-3\geq 0.$ Therefore, when $d=1$, that is, $D$ is a line, we cannot apply our method to construct closed two-forms on the moduli space $\mathcal{M}^*(\mathbb{CP}^2;F_D,r,c\,l_\infty,n)^{sm}.$

Fix $d=2.$ In this case $D$ is a nondegenerate conic in $\mathbb{CP}^2.$ Let $l$ be a line in $\mathbb{CP}^2$ of equation $\omega_l=0$, where $\omega_l\in\mathrm{H}^0(\mathbb{CP}^2,\omega_{\mathbb{CP}^2}(2\, D))=\mathrm{H}^0(\mathbb{CP}^2,\mathcal{O}_{\mathbb{CP}^2}(1)).$ Since for locally free sheaves on $\mathbb{CP}^1$ the triviality condition is equivalent to the semistability condition, we can define the open subset $\mathcal{M}_{lf}^*(\mathbb{CP}^2;F_D,r,c\,l_\infty,n)\subset \mathcal{M}^*(\mathbb{CP}^2;F_D,r,c\,l_\infty,n)$ consisting of isomorphism classes of $(D,F_D)$-framed vector bundles $[(E,\alpha)]$ on $\mathbb{CP}^2$, with $E$ trivial along $l.$

Let $[(E,\alpha)]$ be a point in $\mathcal{M}_{lf}^*(\mathbb{CP}^2;F_D,r,c\,l_\infty,n)^{sm}.$ Consider the induced map on the Ext-group given by the multiplication by $\omega_l$:
\begin{equation*}
 (\omega_l)_*\colon \mathrm{Ext}^1(E,E(-D))\longrightarrow \mathrm{Ext}^1(E,E\otimes \omega_X(D)),
\end{equation*}
that is,
\begin{equation*}
 (\omega_l)_*\colon \mathrm{Ext}^1(E,E(-2))\longrightarrow \mathrm{Ext}^1(E,E(-1)).
\end{equation*}
If we tensor by $E$ and then apply the functor $\mathrm{Hom}(E,\cdot)$ to the short exact sequence
\begin{equation*}
0\longrightarrow \mathcal{O}_{\mathbb{CP}^2}(-2)\stackrel{\cdot \omega_l}{\longrightarrow} \mathcal{O}_{\mathbb{CP}^2}(-1)\longrightarrow \mathcal{O}_l(-1)\longrightarrow 0,
\end{equation*}
we get the long exact sequence
\begin{equation*}
\cdots \rightarrow \mathrm{Hom}(E\vert_l,E\vert_l(-1))\rightarrow \mathrm{Ext}^1(E,E(-2))\stackrel{(\omega_l)_*}{\rightarrow} \mathrm{Ext}^1(E,E(-1))\rightarrow \mathrm{Ext}^1(E\vert_l,E\vert_l(-1))\rightarrow \cdots
\end{equation*}
Since $E$ is trivial along $l$ and $\mathrm{H}^0(l,\mathcal{O}_l(-1))=\mathrm{H}^1(l,\mathcal{O}_l(-1))=0$, we get that $(\omega_l)_*$ is an isomorphism, hence $\tau(\omega_l)_{[(E,\alpha)]}$ is non-degenerate at $[(E,\alpha)]$ by Proposition \ref{prop:non-dege}. Therefore $\tau(\omega_l)$ defines a holomorphic symplectic structure on $\mathcal{M}_{lf}^*(\mathbb{CP}^2;F_D,r,c\,l_\infty,n)^{sm}.$ For $F_D\cong \mathcal{O}_D^{\oplus r}$, we get $c=0$ and $\mathcal{M}_{lf}^*(\mathbb{CP}^2;\mathcal{O}_D^{\oplus r},r,0,n)$ is a holomorphic symplectic variety of dimension $2rn.$
\begin{remark}
In this case $\omega^{-1}_{\mathbb{CP}^2}(-2\,D)\cong \mathcal{O}_{\mathbb{CP}^2}(-1)$, which has not global sections. Hence it is not possible to use Bottacin's result to construct (nondegenerate) Poisson structures on $\mathcal{M}_{lf}^*(\mathbb{CP}^2;\mathcal{O}_D^{\oplus r},r,0,n).$ \triend
\end{remark}

\subsection*{The Hirzebruch surfaces}

Let $p$ be a nonnegative integer number. We denote by $\mathbb{F}_p$ the $p$-th Hirzebruch surface $\mathbb{F}_p:=\mathbb{P}(\mathcal{O}_{\mathbb{CP}^1}\oplus \mathcal{O}_{\mathbb{CP}^1}(-p))$, which is the projective closure of the total space of the line bundle $\mathcal{O}_{\mathbb{CP}^1}(-p)$ on $\mathbb{CP}^1.$ One can describe explicitly $\mathbb{F}_p$ as the divisor in $\mathbb{CP}^2\times \mathbb{CP}^1$
\begin{equation*}
\mathbb{F}_p=\{([z_0:z_1:z_2],[z:w])\in \mathbb{CP}^2\times \mathbb{CP}^1\,\vert\, z_1w^p=z_2 z^p\}.
\end{equation*}
Let us denote by $\pi\colon \mathbb{F}_p \rightarrow \mathbb{CP}^2$ the projection onto $\mathbb{CP}^2$ and by $l_\infty$ the inverse image of a generic line of $\mathbb{CP}^2$ through $\pi.$ Then $l_\infty$ is a smooth connected big and nef curve of genus zero. The Picard group of $\mathbb{F}_p$ is generated by $l_\infty$ and the fibre $F$ of the projection $\mathbb{F}_p\rightarrow \mathbb{CP}^1.$ One has
\begin{equation*}
l_\infty^2=p,\; l_\infty\cdot F=1,\; F^2=0.
\end{equation*}
In particular, the canonical divisor $K_p$ can be expressed as
\begin{equation*}
K_p=-2\,l_\infty +(p-2)\,F.
\end{equation*}

Let us consider the case $p=0.$ In this case $\mathbb{F}_0$ is the ruled surface $\mathbb{CP}^1\times \mathbb{CP}^1.$ Moreover, $K_0=-2\,l_\infty -2F.$ By \cite[Corollary V-2.18]{book:hartshorne1977}, there exists a smooth connected ample curve $D$ in the complete linear system $\vert l_\infty +F\vert.$ By the adjuction formula (see, e.g.,  \cite[Proposition V-1.5]{book:hartshorne1977}), $D$ has genus zero.

Let $n\in \mathbb{Z}$ and $F_D$ a Gieseker-semistable locally free $\mathcal{O}_D$-module of rank $r$ and degree $a+b$, for $a,b\in \mathbb{Z}$, on $D.$ Let $\mathcal{M}^*(\mathbb{F}_0;F_D,r,a\,l_\infty+b\,F,n)$ be the moduli space of $(D,F_D)$-framed sheaves on $\mathbb{F}_0$ of rank $r$, first Chern class $a\,l_\infty+b\,F$ and second Chern class $n.$ Since $K_{\mathbb{F}_0}+2\,D=0$, the line bundle $\omega_{\mathbb{F}_0}(2\,D)$ is trivial and for $1\in \mathrm{H}^0(\mathbb{F}_0,\omega_{\mathbb{F}_0}(2\,D))\cong \mathbb{C}$, the two-form $\tau(1)$ defines a holomorphic symplectic structure on $\mathcal{M}^*(\mathbb{F}_0;F_D,r,a\,l_\infty+b\,F,n)^{sm}.$ If $F_D\cong \mathcal{O}_D^{\oplus r}$, we have $b=-a$ and the moduli space $\mathcal{M}^*(\mathbb{F}_0;F_D,r,a(l_\infty-F),n)$ is a holomorphic symplectic variety of dimension $2(rn+(r-1)a^2).$

Consider $p=1.$ The first Hirzebruch surface $\mathbb{F}_1$ is isomorphic to the blowup of $\mathbb{CP}^2$ at a point. Consider the complete linear system $\vert l_\infty + F\vert.$ Again, by \cite[Corollary V-2.18]{book:hartshorne1977}, there exists a smooth connected curve $D$ on $\vert l_\infty + F\vert.$ By \cite[Theorem V-2.17]{book:hartshorne1977}, $\vert l_\infty + F\vert$ is base-point-free, hence $D$ is nef. Since $D^2=3$, $D$ is also big. By the adjuction formula, the genus of $D$ is zero. Moreover, from $K_1=-2\,l_\infty-F$, it follows that $K_1+2\,D=F.$

Let $n\in \mathbb{Z}$ and $F_D$ a Gieseker-semistable locally free $\mathcal{O}_D$-module of rank $r$ and degree $2a+b$, for $a,b \in \mathbb{Z}$, on $D.$ Let $\mathcal{M}^*(\mathbb{F}_1;F_D,r,a\,l_\infty+b\,F,n)$ be the moduli space of $(D,F_D)$-framed sheaves on $\mathbb{F}_1$ of rank $r$, first Chern class $a\,l_\infty +b\,F$ and second Chern class $n.$ Let $l$ be a smooth connected curve of genus zero in $\mathbb{F}_1$ defined by a nonzero section $\omega_l\in \mathrm{H}^0(\mathbb{F}_1,\omega_{\mathbb{F}_1}(2\, D)).$ As in the previous example, the two-form $\tau(\omega_l)$ defines a holomorphic symplectic structure on the smooth locus of the moduli space $\mathcal{M}_{lf}^*(\mathbb{F}_1;F_D,r,a\,l_\infty+b\,F,n)$ parametrizing isomorphism classes of $(D,F_D)$-framed vector bundles $[(E,\alpha)]$ on $\mathbb{F}_1$, with $E$ trivial along $l.$ For $F_D\cong \mathcal{O}_D^{\oplus r}$, $b=-2a$ and $\mathcal{M}_{lf}^*(\mathbb{F}_1;F_D,r,a\,l_\infty-2a\,F,n)$ is a holomorphic symplectic variety of dimension $2rn+3(r-1)a^2.$
\begin{remark}
As in the previous example, this case is not covered by Bottacin's result since $\omega_{\mathbb{F}_1}^{-1}(-2\, D)\cong \mathcal{O}_{\mathbb{F}_1}(-F)$, which has not global sections. \triend
\end{remark}

Finally, let $p=2.$ In this case $\mathbb{F}_2$  is the projective closure of the cotangent bundle $T^*\mathbb{CP}^1$ of the complex projective line $\mathbb{CP}^1.$ Let $D=l_\infty$, $n\in \mathbb{Z}$ and $F_D$ a Gieseker-semistable locally free $\mathcal{O}_D$-module of rank $r$ and degree $2a+b$, with $a,b \in \mathbb{Z}$, on $D.$ 

Let $\mathcal{M}^*(\mathbb{F}_2;F_D,r,a\,l_\infty+b\,F,n)$ be the moduli space of $(D,F_D)$-framed sheaves on $\mathbb{F}_2$ of rank $r$, first Chern class $al_\infty +bF$ and second Chern class $n.$ The canonical divisor of $\mathbb{F}_2$ is $K_2=-2\,D$, hence the line bundle $\omega_{\mathbb{F}_2}(2\, D)$ is trivial and, for $1\in \mathrm{H}^0(\mathbb{F}_2,\omega_{\mathbb{F}_2}(2\,D))\cong \mathbb{C}$, the two-form $\tau(1)$ defines a symplectic structure on $\mathcal{M}^*(\mathbb{F}_2;F_D,r,a\,l_\infty+b\,F,n)^{sm}.$

If $F_D\cong \mathcal{O}_D^{\oplus r}$, we have $b=-2a.$ Let us define $C=l_\infty-2F.$ This is the only irreducible curve in $\mathbb{F}_2$ with negative self intersection. We can normalize the value $a$ in the range $0\leq a \leq r-1$ upon twisting by $\mathcal{O}_{\mathbb{F}_2}(C).$ By \cite[Theorem 3.4]{phd:rava2012}, the moduli space $\mathcal{M}^*(\mathbb{F}_2;F_D,r,a\,C,n)$ is nonempty if and only if $n+a(a-1)\geq 0$, and if this is the case, $\mathcal{M}^*(\mathbb{F}_2;F_D,r,a\,C,n)$ is a holomorphic symplectic variety of dimension $2(rn+(r-1)a^2).$

\subsection*{A ruled surface over an elliptic curve}

Let $C$ be an elliptic curve. Then $C\times \mathbb{CP}^1$ with its first projection is a ruled surface. By \cite[Proposition V-2.8]{book:hartshorne1977}, there exists a section $C_0$ of $C\times \mathbb{CP}^1$, which is an ample divisor. Let $F$ denote the fibre of the projection $C\times \mathbb{CP}^1\rightarrow C.$ Then $\mathrm{Pic}(C\times \mathbb{CP}^1)=\mathbb{Z}C_0\oplus \mathbb{Z}F$, with $C_0^2=F^2=0$, and $C_0\cdot F=1.$ Moreover, $K_{C\times \mathbb{CP}^1}=-2\,C_0.$ 

Let $D$ be a smooth connected ample curve in the complete linear system $\vert C_0\vert.$ By the adjuction formula, the genus of $D$ is one. Let $n\in \mathbb{Z}$ and $F_D$ a Gieseker-semistable locally free $\mathcal{O}_D$-module of rank $r$ and degree $b\in \mathbb{Z}$ on $D.$ 

Let $\mathcal{M}^*(C\times \mathbb{CP}^1;F_D,r,a\,C_0+b\,F,n)$ be the moduli space of $(D,F_D)$-framed sheaves on $C\times \mathbb{CP}^1$ of rank $r$, first Chern class $a\,C_0 +b\,F$, with $a\in \mathbb{Z}$, and second Chern class $n.$ Also in this case, $K_{C\times \mathbb{CP}^1}+2\,D=0$, hence for $1\in \mathrm{H}^0(C\times \mathbb{CP}^1, \omega_{C\times \mathbb{CP}^1}(2\,D))\cong \mathbb{C}$, the two-form $\tau(1)$ defines a holomorphic symplectic structure on $\mathcal{M}^*(C\times \mathbb{CP}^1;F_D,r,a\,C_0+b\,F,n)^{sm}.$

\end{document}